\documentclass[12pt, english]{article}
\RequirePackage[OT1]{fontenc}
\RequirePackage{amsthm,amsmath,natbib}
\usepackage{calc, xspace}
\usepackage{ifpdf}
\usepackage{graphicx}
\usepackage[margin=1in]{geometry}
\RequirePackage{amsfonts, euscript}
\RequirePackage[colorlinks,citecolor=blue,urlcolor=blue]{hyperref}
\RequirePackage{framed}
\RequirePackage[small]{caption}
\usepackage[utf8x]{inputenc}
\usepackage{dsfont}
\usepackage{bm}
\newcommand{\Reals}{\mathbb{R}} 
\usepackage{url}
\usepackage{comment}
\usepackage[british]{babel}
\usepackage{setspace}
\usepackage{sectsty}
\usepackage{titlesec}
\titleformat{\section}[hang]{\centering \bfseries\large}{\thesection.}{0.4em}{}

\usepackage{float}
\usepackage{color,soul}

 \newcommand{\eps}{\varepsilon}
\newcommand{\Prob}[1]{\operatorname{\mathbb{P}}\left[#1\right]}

 \newcommand{\Expect}[1]{\operatorname{\mathbb{E}}\left[#1\right]}
\newcommand{\parens}[1]{\left(#1\right)}
\newcommand{\beqn}{\vspace{-0.25cm}\begin{eqnarray*}}
\newcommand{\eeqn}{\end{eqnarray*}}
\newcommand{\bneqn}{\vspace{-0.25cm}\begin{eqnarray}}
\newcommand{\eneqn}{\end{eqnarray}}
\newcommand{\bracks}[1]{\left[#1\right]}
\newcommand{\expesub}[2]{\mathbb{E}_{#1}\bracks{#2}}
\newcommand{\bs}{\textbf{s}}
\newcommand{\bt}{\textbf{t}}
\newcommand{\MEXIT}{\emph{MEXIT}\xspace}
\newcommand{\cadlag}{c\`adl\`ag\xspace}

\newcommand{\Diagonal}[1]{\Delta_{#1}}
\newcommand{\DiagonalInjection}[1]{\operatorname{\textnormal{\i}}_{\Diagonal{#1}}} %%% QZ changed \text{\i} to \mathrm{\i} so that the symbol looks the same in text and in equation. %%% WSK changed to \textnormal: \mathrm creates errors!  
\newcommand{\SigmaAlgebra}{\mathcal{E}}
\renewcommand{\P}{\mathbb{P}}
\newcommand{\StateSpace}{E}
\newcommand{\TV}{\textsc{TV}}
\def\tcr{\textcolor{black}}

\def\bs{{\bf s}}

\def\alphah{\widehat{\alpha}}
\def\alphat{\widetilde{\alpha}}
\def\Pt{\widetilde{P}}
\def\Normal{{\rm Normal}}
 \newcommand{\Law}[1]{\mathcal{L}\left({#1}\right)}
 \newcommand{\Leb}{\operatorname{Leb}}
 
  \renewcommand{\d}{\,\operatorname{d}}
   \DeclareMathOperator{\dist}{dist}

\theoremstyle{plain}
 \newtheorem{thm}{Theorem}
 \newtheorem{proposition}[thm]{Proposition}
 \newtheorem{lem}[thm]{Lemma}
 \newtheorem{cor}[thm]{Corollary}
 
\theoremstyle{definition}
 \newtheorem{defn}[thm]{Definition}
\theoremstyle{remark}
 \newtheorem{rem}[thm]{Remark}

\title{MEXIT: Maximal un-coupling times for stochastic processes}
\author{%
 Philip A. Ernst\footnote{\href{mailto:philip.ernst@rice.edu}{\texttt{philip.ernst@rice.edu}}},
 Wilfrid S.~Kendall\footnote{\href{mailto:w.s.kendall@warwick.ac.uk}{\texttt{w.s.kendall@warwick.ac.uk}}},
 Gareth O.~Roberts\footnote{\href{mailto:gareth.o.roberts@warwick.ac.uk}{\texttt{gareth.o.roberts@warwick.ac.uk}}},
 Jeffrey S.~Rosenthal\footnote{\href{mailto:jeff@math.toronto.edu}{\texttt{jeff@math.toronto.edu}}}  
} 
\begin{document}
 \maketitle

%%%%%%%%%%%%%%%%%%%%%%%%%%%%%%%%
%% Title page

% =========================================   
\begin{abstract}
% \tcr{
Classical coupling constructions arrange for copies of the \emph{same} Markov process started at two \emph{different}
initial states to become equal as soon as possible. In this paper, we consider an alternative coupling framework
in which one seeks to arrange for two \emph{different} Markov (or other stochastic) processes to remain equal
for as long as possible, when started in the \emph{same} state. We refer to this ``un-coupling'' or ``maximal agreement'' construction as \emph{MEXIT}, standing for ``maximal exit''. After highlighting the importance of un-coupling arguments in a few key statistical and probabilistic settings, we develop an explicit \MEXIT construction for
stochastic processes in discrete time with countable state-space. This construction is generalized to random processes on general state-space running in continuous time, and then exemplified by discussion of \MEXIT for Brownian motions with two different constant drifts.\\
\end{abstract}
\textsc{Fr\'echet class;}
\textsc{Hahn-Jordan decomposition;}
\textsc{Markovian coupling;} 
\textsc{MCMC;}  
\textsc{meet measure;}
\textsc{MEXIT;}
\textsc{one-step Minorization;}
\textsc{pseudo-marginal MCMC;}
\textsc{recognition lemma for maximal coupling;}
\textsc{un-coupling.}
% =========================================   

 \section{Introduction}\label{sec:introduction}

Coupling is a 
device commonly employed 
in probability theory for learning about
distributions of certain random variables by means of judicious construction in ways which depend on other random variables (\citet{Lindvall-1992} and \citet{Thorisson-2000}). 
Such coupling constructions are often used to prove convergence of
Markov processes to stationary
distributions (\citet{Pitman-1976}), especially for Markov chain
Monte Carlo (MCMC) algorithms (\citet[and references therein]{RobertsRosenthal-2004}), by seeking to build two different
copies of the \emph{same} Markov process started at two \emph{different}
initial states in such a way that they become equal at a fast rate.
Fastest possible rates are achieved by the \emph{maximal coupling} constructions which were introduced and studied in \citet{Griffeath-1975}, \citet{Pitman-1976}, and \citet{Goldstein-1978}.
Our results and methods are closely related to 
the work of \citet{Goldstein-1978}, who deals with rather general discrete-time random processes.
Our situation is related to a time-reversal of the situation studied by \citet{Goldstein-1978}. 
However
% we are able to avoid some measure-theoretic assumptions, since we work initially with finite time intervals and focus on explicit construction of probability measures.\\
our approach seems more direct.
\\
%Additionally we are able to make a clear link with copula theory.
%\cite{Pitman-1976} provided a simplified construction of Griffeath's maximal coupling, offering an explicit construction for discrete-state-space discrete-time Markov chains. 
\indent In the current work, we consider what might be viewed as the dual problem where coupling is
used to try to construct two \emph{different} Markov (or
other stochastic) processes which remain equal
for as long as possible, when they are started in the \emph{same} state.
That is, 
we move from consideration of the coupling time to focus on
% instead we consider
the \emph{un-coupling time} at which the processes
diverge, and try to make that as \emph{large} as possible. 
We refer to this as \emph{MEXIT} (standing for ``maximal exit'' time). 
While finalizing our current work, it came to our attention that this construction is the same as the \emph{maximal agreement coupling time} of the August 2016 work of \citet{Vollering-2016},
who additionally derives a lower bound on \MEXIT.
% \tcr{
Nonetheless, we believe the current work complements \citet{Vollering-2016} well. 
% who additionally derives a lower bound on \MEXIT. 
It offers an explicit treatment of discrete-time countable-state-space, generalizes the continuous-time case, 
and discusses a number of significant applications of \MEXIT. We note that the work of \citet{Vollering-2016} does not consider the continuous-time case. 
% and provides a detailed discussion of \MEXIT for general random processes and for diffusions.
% }

In addition to being a natural mathematical question, \MEXIT has direct applications 
to many key statistical and probabilistic settings (see Section \ref{sec:applic} below). 
In particular, couplings which are \emph{Markovian} and \emph{faithful}
(\citet{Rosenthal-1997}, i.e.\ couplings which preserve the marginal
update distributions even when conditioning on both processes;
alternatively ``co-adapted'' or ``immersion'', 
depending on the extent to which one wishes to emphasize the underlying filtration as in \citet{BurdzyKendall-2000} and \citet{Kendall-2013a}) 
are the most
straightforward to construct, but often are \emph{not} maximal, % and
while
more complicated non-Markovian and non-faithful couplings lead to
stronger bounds.
The same is true in the context of \MEXIT. 
%    \TBC{
% %    Motivate by notion of ``efficient coupling'' \cite{BurdzyKendall-2000,BanerjeeKendall-2015}.
% %    Note its relationship to the Aldous coupling inequality.
%    Note relationship of successful coupling to the Aldous coupling inequality.
% %    \\
% %    Also mention faithful / co-adapted / Markovian / immersion couplings \cite{Rosenthal-1997,BurdzyKendall-2000,Kendall-2013a}.
% %    (Chance to set down in print that these are all more-or-less the same idea!)
%    \\
%    Notion of ``\MEXIT'' (maximal exit). 
%    Is \MEXIT always well-defined? (YES!) And what about ``multiple \MEXIT''?
%    \\
%    Refer generously to \cite{Vollering-2016}, 
%    note the terminology of ``maximal agreement coupling'', also the lower bound on \MEXIT time,
%    and make it clear how our work goes further than that of \citeauthor{Vollering-2016}.
%    }
% %    \\
% %    Include and expand on introductory material from Ernst-Roberts-Rosenthal, June 2016.
% %    Especially, develop Section 2 of notes by Ernst et al., on Applications. (Or separate section?)

% =========================================   
\section{Applications}\label{sec:applic}
To motivate the natural role of \MEXIT in the existing literature, we first consider the role of un-coupling arguments in a few statistical and probabilistic settings.

% =========================================   
\subsection{Bounds on accuracy for statistical tests}\label{sec:testing}

Un-coupling 
% plays a prominent role in 
has an impact on the theory
classical statistical testing. In \citet{Farrell-1964}, amongst other sources, some function of the data (but not the data itself) is assumed to have been observed.  
% \tcr{
A statistical test is then constructed to enable detection of the distribution from which
the observed data have been sampled.
% }
For example, suppose that data $X_1,X_2,\ldots$ are generated
as a draw either from a multivariate probability distribution $\P_1$ or from a multivariate probability distribution $\P_2$.  
% \tcr{
The goal is to determine whether the data was generated from $\P_1$ or from $\P_2$. For some function $h$ of the data, and some acceptance region
$A$, the statistical test 
decides in favor of $\P_1$ if $h(X_1,\ldots,X_n) \in A$ and otherwise decides in favor of
$\P_2$.
% } 
% Further suppose that we have the existence of

Suppose that there exists
an un-coupling time $T$, such that if
$X_1,X_2,\ldots$ are generated from $\P_1$, and if $Y_1,Y_2,\ldots$
are generated from $\P_2$
% . Let
then it is exactly the case that
$X_i=Y_i$ for all $1 \le i \le T$
(so that \(X_i\neq Y_i\) for all \(i > T\)).
We use \(\P\) to refer to the joint distribution (in fact, the coupling) of $\P_1$ and $\P_2$.

The following proposition uses the un-coupling probabilities to %provide
recover 
a lower bound on the accuracy of such statistical tests
related to \citet[Theorem 1]{Farrell-1964}.
\begin{proposition} Under the above assumptions,
the sum of the probabilities of Type-I and Type-II errors of our statistical
test is at least $\Prob{T>n}$.
\end{proposition}

\begin{proof}
We apply elementary arguments to
the sum of the probabilities of Type-I and Type-II errors:
\begin{align*}
& \P_2[h(Y_1,\ldots,Y_n) \in A]
+ \P_1[h(X_1,\ldots,X_n) \not\in A]
\quad=\quad
% $$
% $$
\\
&\quad = \quad
\P_2[h(Y_1,\ldots,Y_n) \in A]
+ 1 - \P_1[h(X_1,\ldots,X_n) \in A]
% $$
% $$
\\
&\quad = \quad
1 - \Big(\P_{1}[h(Y_1,\ldots,Y_n) \in A]
- \P_{2}[h(X_1,\ldots,X_n) \in A]\Big)
% $$
% $$
\\
&\quad\ge\quad
1 - \left|\P_{1}[h(Y_1,\ldots,Y_n] \in A]
- \P_{2}[h(X_1,\ldots,X_n) \in A]\right|
% \quad \text{ \tcr{(triangle inequality) (PE1: remove? this is merely $x\leq |x|$ )}}
% $$
% $$
\\
&\quad\ge\quad
1 - \|\mathcal{L}_{\P_1}(X_1,\ldots,X_n) - \mathcal{L}_{\P_2}(Y_1,\ldots,Y_n)\|
\quad \text{ (definition of total variation distance)}
% $$
% $$
\\
&\quad\ge\quad
1 - \Prob{X_i\not=Y_i \ {\rm for \ some} \ 1 \le i \le n}
\quad \text{ (coupling inequality)}
% $$
% $$
\\
&\quad = \quad
1 - (1-\Prob{X_i=Y_i \text{ for all } \ 1 \le i \le n})
% $$
% $$
\\
&\quad = \quad
\Prob{X_i=Y_i \ {\rm for \ all} \ 1 \le i \le n}
\quad = \quad
\Prob{T>n}
\, .
% $$
% $$
% \ \ge \
% 1 - [1-\P(T>n)]
% \ = \
% \P(T>n)
% \, .
\end{align*}
\end{proof}
% \tcr{
% The above leads to a powerful conclusion regarding the accuracy of statistical tests. 
% If we know that $\Prob{T>t}=1-\epsilon$, \textit{then there is no test} which can tell the two processes apart with accuracy greater than $\epsilon$!
% }

\subsection{Two independent coin flips}\label{sec:flips}
% \tcr{
We now turn to the classical probabilistic  paradigm of coin flips.
% }. 
Let 
\(X\) and \(Y\) represent
% $\{X_n\}$ and $\{Y_n\}$ to be 
two different
sequences of i.i.d.\ coin flips, with probabilities of landing on H (heads) to be $q$ and
$r$ respectively, where $0 \leq r \leq q \leq 1/2$.
Suppose that we wish to maximise the 
% variable
length of the initial segment for which coin flips agree:
$$
S \quad = \quad \max\{t : X_i=Y_i \ \hbox{\rm for all} \ 1 \le i \le t\}
\, .
$$
% \tcr{Throughout this section, for purposes of concreteness}, we let 
For concreteness, we will set
$q=0.5$ and
$r=0.4$
throughout this section; the generalization to other values is immediate.

% =========================================   
\subsubsection{Markovian Faithful Coupling for Independent Coin Flips}\label{sec:flips-faithful}
The ``greedy'' (Markovian and faithful) coupling carries out the best
``one-step minorization'' coupling possible, separately at each iteration.
One-step minorization is essentially maximal coupling for single random variables.
In this case, that means that for each flip, $\Prob{X=Y=H} = 0.4$,
$\Prob{X=Y=T} = 0.5$, and $\Prob{X=H, \ Y=T} = 0.1$.  
This preserves the
marginal distributions of $X$ and $Y$, and yields $\Prob{X=Y}=0.9$ at each
step.  
Accordingly, the probability of agreement continuing for at least $n$ steps is given by $\Prob{X_i=Y_i \text{ for } 1 \le i \le n}
= (0.9)^n$.

% =========================================   
\subsubsection{A Look-ahead Coupling for Independent Coin Flips}\label{sec:flips-look-ahead1}
Let a ``look-ahead'' coupling be a coupling which instead uses an \(n\)-step minorization
couple on the entire sequence of \(n\) coin tosses, so that for each sequence
$s$ of $n$ different Heads and Tails, it sets $\Prob{X=Y=s} = \min(\Prob{X=s}, \,
\Prob{Y=s})$.  Consequently, if $s$ has $m$ Heads and $n-m$ Tails, then
\[
 \Prob{X=Y=s} \quad=\quad \min\{0.5^n, \, 0.4^m 0.6^{n-m}\}\,.
\]
Elementary events for which \(X\) and \(Y\) disagree are assigned probabilities which
preserve the marginal distributions of $X$ and of $Y$. The simplest
way to implement this is to use ``independent residuals'', but other choices
are also possible.
% , corresponding to ``discrete copulas''.

This look-ahead coupling leads to a
larger probability that $X=Y$.  Indeed, even in the case $n=2$,
the probability of agreement over two coin flips under the greedy coupling is given by 
\begin{equation*}
\Prob{X=Y}\quad= \quad (0.9)^2 = 0.81\,.
\end{equation*}
The look-ahead coupling delivers a strictly larger probability of agreement over two coin flips:
\begin{align*}
\Prob{X=Y} &\quad=\quad \min(0.5^2, 0.4^2) + \min(0.5^2, 0.6^2)
+ 2 \min(0.5^2, 0.4 \cdot 0.6)
\\
&\quad=\quad 0.4^2 + 0.5^2 + 2 \cdot 0.4 \cdot 0.6
= 0.16 + 0.25 + 0.48=0.89
\,.
\end{align*}
% which is larger. 
When $n=2$, the matrix of joint probabilities for $X$ and $Y$
under the look-ahead coupling is calculated to be:
$$
\begin{tabular}{|c||c|c|c|c||c|}
\hline
$X \backslash Y$ & HH & HT & TH & TT & SUM \cr
\hline
HH & 0.16 & 0 & 0 & 0.09 & 0.25 \cr
HT & 0 & 0.24 & 0 & 0.01 & 0.25 \cr
TH & 0 & 0 & 0.24 & 0.01 & 0.25 \cr
TT & 0 & 0 & 0 & 0.25 & 0.25 \cr
\hline
SUM & 0.16 & 0.24 & 0.24 & 0.36 & 1 \cr
\hline
\end{tabular}
$$

Marginalizing this coupling on the initial coin flip (``projecting back'' to the initial flip, with $n=1$), we see
that $\Prob{X_1=Y_1=H} = 0.16 + 0.24 = 0.4$, and
$\Prob{X_1=Y_1=T} = 0.24 + 0.01 + 0.25 = 0.5$,
and $\Prob{X_1=H, \ Y_1=T} = 0.09 + 0.01 = 0.1$. 
The projection to
the initial flip 
yields the same agreement probability as would have been attained by maximizing
the probability of staying together for just one flip ($n=1$).  
That is, the
$n=2$ look-ahead coupling construction is \emph{compatible} with the $n=1$
construction.  
% [More about this later!]

Finally, it is worth noting that the $n=2$ look-ahead coupling is certainly not faithful.  For
example, $\Prob{X_2=H \, | \, X_1=Y_1=H} = 0.4 \neq 0.5$,
and $\Prob{X_2=H \, | \, X_1=H, \ Y_1=T} = 0.9 \neq 0.5$, \emph{etc}.

% =========================================   
\subsection{A Look-ahead coupling for independent coin flips: the case $n=3$}\label{sec:flips-look-ahead2}
\def\sp{\underline{ \ \ }}
%
% When $n=3$, t
The matrix of joint probabilities for $X$ and $Y$
under the look-ahead coupling for \(n=3\) is more complicated, but
can be calculated as:
$$
\begin{tabular}{|c||c|c|c|c|c|c|c|c||c|}
\hline
$X \backslash Y$ & HHH & HHT & HTH & HTT & THH & THT & TTH & TTT & SUM \cr
\hline
HHH & 0.064 & 0 & 0 & 0.0078 & 0 & 0.0078 & 0.0078 & 0.0375 & 0.125 \cr
HHT & 0 & 0.096 & 0 & 0.0037 & 0 & 0.0037 & 0.0037 & 0.0178 & 0.125 \cr
HTH & 0 & 0 & 0.096 & 0.0037& 0 & 0.0037 & 0.0037 & 0.0178 & 0.125 \cr
HTT & 0 & 0 & 0 & 0.125 & 0 & 0 & 0 & 0 & 0.125 \cr
THH & 0 & 0 & 0 & 0.0037 & 0.096 & 0.0037 & 0.0037 & 0.0178 & 0.125 \cr
THT & 0 & 0 & 0 & 0 & 0 & 0.125 & 0 & 0 & 0.125 \cr
TTH & 0 & 0 & 0 & 0 & 0 & 0 & 0.125 & 0 & 0.125 \cr
TTT & 0 & 0 & 0 & 0 & 0 & 0 & 0 & 0.125 & 0.125 \cr
\hline
SUM & 0.064 & 0.096 & 0.096 & 0.144 & 0.096 & 0.144 & 0.144 & 0.216 & 1 \cr
\hline
\end{tabular}
$$
 With these probabilities, we compute that
$$
\Prob{X=Y}= 0.064 + 3 \times 0.096 + 4 \times 0.125
= 0.852
\, .
$$
This is greater than the agreement probability of $0.9^3 = 0.729$ that would have be achieved
via the greedy coupling.
% (More generally, we think that perhaps
% it is always possible to ensure compatibility in this sense, but we
% are not yet sure about that.)
It is natural to wonder whether or not it is possible always to ensure that such a construction works not just for one fixed time but for all times.
We further expound on this point in Sections \ref{sec:MEXIT} and \ref{sec:generalMEXIT}, where discussion of a much more general context shows that that such constructions always exist.

% =========================================   
\subsubsection{Optimal Expectation}\label{sec:optimal-expect}
Until now, this section has focused on maximising $\Prob{X_i=Y_i \ \text{ for all }
\ 1 \le i \le n}$, 
which is to say, maximizing $\Prob{S \ge n}$ with $S$ being the time of first disagreement as above. 
We now consider the related question of maximizing the expected value $\Expect{S}$ . 
Using the greedy coupling, clearly
$$
\Expect{S} \quad=\quad \sum_{j=1}^\infty \Prob{S \ge j}
\quad=\quad \sum_{j=1}^\infty 0.9^j
\quad=\quad 0.9/(1-0.9) = 9
\, .
$$
If the different look-ahead couplings are chosen to be compatible, then
% we can use a similar sum to show
this shows
that $\Expect{S}$ is the sum  for \(r=1,2,\ldots\) of the
probabilities that the $j^{\rm th}$ look-ahead coupling was successful.
% This appears(?) to be what was used on page~3 of the Ernst/Shepp notes.
% However, it does require justification, which appears to be lacking.
% 
% OPEN PROBLEM: Can the above look-ahead couplings always be made to be
% compatible in this sense?
% (If yes, then we get the above improved bound on $\Expect{S}$, but if not
% then not.)
The work of Sections \ref{sec:MEXIT} and \ref{sec:generalMEXIT} shows that such a choice is always feasible, even for very general random processes indeed.

% =========================================   
\subsection{Adaptive MCMC}\label{sec:adaptive}
% We highlight the use of u
Un-coupling arguments play a natural role in the adaptive MCMC (Markov-chain Monte Carlo) literature, highlighted in particular by the work of \citet{RobertsRosenthal-2007}. \citet{RobertsRosenthal-2007} prove
convergence of \emph{adaptive} MCMC by comparing an adaptive process
to a process which ``\emph{stops} adapting'' at some point, and then by showing that the
two processes have a high probability of remaining equal long enough
such that the second process (and hence also the first process) converge
to stationarity. The authors accomplish this by considering a sequence of adaptive Markov
kernels $P_{\Gamma_1},P_{\Gamma_2},\ldots$ on a state space $\mathcal{X}$, where
$\{P_\gamma : \gamma\in\mathcal{Y}\}$ are a collection of Markov kernels each
having the same stationary probability distribution $\pi$, and the
$\Gamma_i$ are $\mathcal{Y}$-valued random variables which are ``adaptive'' (i.e.,
they depend on the previous Markov chain values but not on future values).
Under appropriate assumptions, the authors prove that a Markov
chain 
\(X\)
% $\{X_n\}$ 
which evolves \emph{via} the adaptive Markov kernels will still
converge to the specified stationary distribution $\pi$.  The key step in
the proof of the central result \citet[Theorem 5]{RobertsRosenthal-2007} is an un-coupling approach, highlighted below.

% In the proof of
\citet[Theorem 5]{RobertsRosenthal-2007} 
% it is assumed 
assume
that, for any $\eps>0$, there
is a non-negative integer \(N=N(\eps)\)
% $N = N(\eps) \in\Numbers$ 
such that
$$
\|P_\gamma^N(x,\cdot) - \pi(\cdot)\|_\TV
\quad \le \quad \eps
$$
for all $x\in\mathcal{X}$ and $\gamma\in\mathcal{Y}$ (where $\|\cdot\|_\TV$ denotes
total variation norm of a signed measure). Furthermore,
there is 
a non-negative integer \(n^*=n^*(\eps)\)
% $n^* = n^*(\eps) \in \Numbers$ 
such that
with probability at least $1-\eps/N$,
$$
\sup_{x\in\mathcal{X}} \|P_{\Gamma_{n+1}}(x,\cdot)-P_{\Gamma_n(x,\cdot)}\|_\TV \quad \le \quad
\eps/N^2
$$
for all $n \ge n^*$.

These assumptions are used to prove, for any $K \ge n^*+N$,
the existence of a pair of processes
\(X\) and \(X'\)
% $\{X_n\}$ and $\{X'_n\}$ 
defined
for $K-N \le n \le K$, such that 
\(X\)
% $\{X_n\}$ 
evolves \emph{via} the adaptive
transition kernels $P_{\Gamma_n}$, while
\(X'\)
% $\{X'_n\}$ 
evolves \emph{via} the fixed
kernel $P' = P_{\Gamma_{K-N}}$. With probability at least
$1-2\eps$, the two processes remain equal for all times $n$
with $K-N \le n \le K$. Hence, their un-coupling probability over this
time interval is bounded above by $2\eps$. Consequently, conditional on
$X_{K-N}$ and $\Gamma_{K-N}$, the law of $X_K$ lies within $2\eps$
(measured in total variation distance) of the law of $X'_K$, which in turn lies
within $\eps$ of the stationary distribution $\pi$.  Hence,
the law of $X_K$ is within $3\eps$ of $\pi$.  Since
this holds for any $\eps>0$ (for sufficiently large
$K=K(\eps)$), it follows that the law of $X_K$ converges
to $\pi$ as $K\to\infty$.  Accordingly the adaptive process
\(X\)
% $\{X_n\}$ 
is indeed a ``valid'' Monte Carlo algorithm for approximately
sampling from $\pi$; namely it converges asymptotically to $\pi$. 
The proof of
a more general result (\citet[Theorem 13]{RobertsRosenthal-2007}), is quite similar, only requiring one
additional~$\eps$.

 \section{\MEXIT for discrete-time countable state-space}\label{sec:MEXIT}

    Having motivated the prominence of un-coupling arguments in key statistical and probabilistic settings, we now turn to an explicit construction of \MEXIT. We begin by considering
    two discrete-time stochastic processes defined on the same countable discrete state-space, begun at the same initial state \(s_0\).
    We extend the state-space by keeping track of the past trajectory of each stochastic process (its ``genealogy'').
    % We suppose that these chains are governed by transition probability kernels \(p(a,b)\) and \(q(a,b)\), respectively.
    % We extend the state-space by keeping track of the past trajectory of each chain (its ``genealogy'').
    The state of one of these stochastic processes at time \(n\) will thus  be a sequence or genealogy \(\bs =(s_0,s_1,\ldots,s_n)\) of \(n+1\) states,
    and these stochastic processes are then time-inhomogeneous Markov chains governed at time \(n\) by transition probability kernels \(p(\bs,b)\) and \(q(\bs,b)\), respectively..
     Let \(\bs\cdot a\) denote the sequence or genealogy \(\bs=(s_0,s_1,\ldots,s_n,a)\) of \(n+2\) states, corresponding to the chain moving to state \(a\) at time \(n+1\).
     Note that if the original processes were originally Markov chains then this notation is equivalent to working with path probabilities \(p(\bs)=p(s_0,s_1)p(s_1,s_2)\ldots p(s_{n-1},s_n)\), \(q(\bs)=q(t_0,t_1)q(t_1,t_2)\ldots q(t_{n-1},t_n)\), with \(p(\bs\cdot a)=p(\bs)p(s_n,a)\) \emph{et cetera}.

     We define a coupling between the two processes as a random process on the Cartesian product of the (extended) state-space with itself,
     whose marginal distributions are those of the individual processes.
    %Since we keep track of the entire genealogies, our approach applies to
    %general discrete-time stochastic processes (not necessarily Markovian).

 \begin{defn}[Coupling of two discrete-time stochastic processes]\label{def:coupling-genealogical}
 A \emph{coupling} of two discrete-time stochastic processes
on a countable state space
with genealogical probabilities $p(\bs)$ and $q(\bt)$ respectively,
 is a random process (\emph{not} necessarily Markov) with state
\((\bs,\bt)\) at time \(n\) given by a pair of genealogies $\bs$
and $\bt$ each of length \(n\),
 such that if the probability of seeing state \((\bs,\bt)\) at time \(n\)
is equal to \(r(\bs,\bt)\), then
  \begin{align}
 \sum_\bt r(\bs, \bt) \quad&=\quad p(\bs) \qquad\text{(row-marginals)}\,,  \label{eq:coup1}  \\   
 \sum_\bs r(\bs, \bt) \quad&=\quad q(\bt) \qquad\text{(column-marginals)}\,.
 \end{align}
Moreover, probabilities at consecutive times are related by
 \begin{equation}\label{eq:inherit}
  \sum_a\sum_b r(\bs\cdot a,\bt\cdot b) \quad=\quad r(\bs,\bt) \qquad\text{(inheritance)}\,.
 \end{equation}
 \end{defn}
 \begin{rem}
  A coupling of two non-genealogical Markov chains can be converted into the above form simply by keeping track of the genealogies.%
% \xNB[WSK1]{I believe \emph{any} random process converts into a (time-inhomogeneous) Markov chain if one keeps track of the genealogy. Thus our results are surprisingly more general than they might appear! (Does this observation lead to an understanding of the original \cite{Griffeath-1975} approach? and, perhaps this is an efficient way to expound the general random process result?)}
 \end{rem}
 \begin{rem}
  We assume that
both processes begin at the same fixed starting point \(s_0\),
so \(p((s_0))=q((s_0))=1\), and
the processes initially have the same trajectory. 
  \MEXIT occurs when first the trajectories split apart and disagree:
the tree-like nature of genealogical state-space means the genealogical
processes will never recombine.
%   \NB[WSK1]{``the tree-like nature of genealogical state-space means they can never then recombine'': but note in conclusion we speculate about immediately trying to re-couple -- this would mean forgetting about the genealogy so this should be omitted or carefully revised. TAKE CARE TO CHECK THIS!}
 \end{rem}
%\begin{rem}
%Astute readers will already be asking whether, in the light of the genealogical nature of Definition \ref{def:coupling-genealogical} the processes really need to be Markov chains.  
%In fact they do not, but we defer consideration of this point to Section \ref{sec:MEXIT}, where we generalize the treatment to apply to rather general random processes both in discrete and in continuous time.
%\end{rem}

 A \MEXIT coupling is one which achieves the bound prescribed by the \citet{Aldous} coupling inequality (Lemma 3.6 therein), thus (stochastically) maximising the time at which the chains split apart.
 \begin{defn}[\MEXIT coupling]\label{def:MEXIT}
  Suppose that the following equation holds for all genealogical states \(\bs\): 
  \begin{equation}\label{eq:MEXIT}
   r(\bs,\bs) \quad=\quad p(\bs)\wedge q(\bs)\,.
  \end{equation}
 Then the coupling is a \emph{maximal exit coupling} (\MEXIT coupling).
 \end{defn}

 We now prove that \MEXIT couplings always exist.
 \begin{thm}\label{thm:MEXIT}
 Consider two discrete-time stochastic processes
taking values in a given countable state-space and started at the same initial state \(s_0\). 
A \MEXIT coupling can always be constructed such that the joint probability $r(\cdot, \cdot)$ satisfies the properties~\eqref{eq:coup1}--\eqref{eq:MEXIT}.   
 \end{thm}
 
\begin{proof}
We claim a MEXIT coupling is given by the following recursive definition
\begin{align*}
  r(\bs \cdot a, \bt\cdot b) = \left\{ \begin{array}{cc}
  p ( \bs \cdot a)  \wedge q (\bs \cdot a)    &   \text{ if } \bt = \bs,  a = b ,  \\
      \pi_1 (\bs \cdot b)      \pi_2 (\bs \cdot a)     \sum_c   d^-(\bs \cdot c)     &   \text{ if } \bt = \bs,  p(\bs ) \leq  q (\bs),  \\ 
      \pi_1 (\bs \cdot b)      \pi_2 (\bs \cdot a)     \sum_c   d^+(\bs \cdot c)     &   \text{ if } \bt = \bs,  p(\bs )  >  q (\bs),  \\  
  \pi_1 (\bt \cdot b)  \pi_2 (\bs \cdot a)  r (\bs, \bt)  &  \text{ if } \bt \neq \bs, 
  \end{array}
  \right. 
\end{align*}
where 
\begin{align*} 
    d^+ ( \bs ) =    ( q(\bs  )   - p (\bs  ) ) \vee 0 ,   &\quad  d^- ( \bs ) =    ( p(\bs  )   - q (\bs  ) ) \vee 0  \\ 
  \pi_1 (\bt \cdot b)  =   \dfrac{  d^+(\bt \cdot b)     }{ \sum_c     d^+(\bt \cdot c)    }   , & \quad 
    \pi_2 (\bs \cdot a)  =    \dfrac{   d^-(\bs \cdot a)     }{ \sum_c  d^-(\bs \cdot c)     }  . 
\end{align*}
We set $\pi_1$ (or $\pi_2$) to zero if the denominator appearing in the definition is zero.
The initial joint probability is given by $r (s_0, s_0) = 1$, which clearly satisfies~\eqref{eq:coup1}--\eqref{eq:MEXIT}.

Now we verify by induction this construction actually satisfies~\eqref{eq:coup1}--\eqref{eq:MEXIT} at each time $n$. 
First, the MEXIT equation~\eqref{eq:MEXIT} holds by construction. 
Second, if $\bs \neq \bt$, we immediately have 
\begin{align*}
 \sum\limits_a \sum\limits_b r(\bs \cdot a, \bt \cdot b) 
=  r(\bs, \bt)
\end{align*}
since  $\sum_c \pi_1 (\bt \cdot c) = \sum_c \pi_2 (\bs \cdot c) = 1$.   
Observe that  
\begin{align*}
 \sum_c   d^-(\bs \cdot c)  +  \sum_c \left(  p ( \bs\cdot c)  \wedge q (\bs \cdot c) \right) 
 =  \sum_c   p ( \bs\cdot c) =  p(\bs), 
\end{align*}
and $d^+(\bs \cdot a) d^-(\bs \cdot a) = 0$. Hence if $p(\bs ) \leq  q (\bs)$, 
\begin{align*}
\sum\limits_a \sum\limits_b r(\bs \cdot a, \bs \cdot b) 
=  \;&  \sum\limits_{c}  \left(  p ( \bs\cdot c)  \wedge q (\bs \cdot c) \right) 
+  \left( \sum_c   d^-(\bs \cdot c)  \right) \sum\limits_{b}\sum\limits_{a \neq b}   \pi_1 (\bs \cdot b)      \pi_2 (\bs \cdot a)     \\
=  \; &  \sum\limits_{c}  \left(  p ( \bs\cdot c)  \wedge q (\bs \cdot c) \right) 
+ \left( \sum_c   d^-(\bs \cdot c)  \right) \dfrac{ \sum\limits_a\sum\limits_b  d^+(\bs \cdot a)  d^-(\bs \cdot b)   }{\left( \sum_c     d^+(\bs \cdot c) \right) \left( \sum_c     d^-(\bs \cdot c) \right) }    \\
= \; &  \sum\limits_{c}  \left(  p ( \bs\cdot c)  \wedge q (\bs \cdot c) \right) 
+ \sum_c   d^-(\bs \cdot c) = p(\bs). 
\end{align*}
Similarly, if $p(\bs) > q(\bs)$, 
\begin{align*}
\sum\limits_a \sum\limits_b r(\bs \cdot a, \bs \cdot b) 
=  q(\bs). 
\end{align*}
Thus we conclude the inheritance property~\eqref{eq:inherit} holds. 
Intuitively, given $r(\bs, \bt)$ at time $n$, we can proceed to time $n+1$ by first filling in the diagonals according to~\eqref{eq:MEXIT}; 
then for each big cell $(\bs, \bt)$, the sum of $ r(\bs \cdot a,\bt \cdot b) $ must be equal to $r(\bs, \bt)$ by~\eqref{eq:inherit} and we fill in all the remaining cells proportionally  by $\pi_1$ and $\pi_2$. 

Now it remains to check the row/column marginal conditions.
We shall only check that the row marginal condition holds. 
If $p(\bs) \leq  q(\bs)$, by the induction assumption, we have $r(\bs, \bs) = p(\bs)$ and $r(\bs, \bt) = $ for any $\bt \neq \bs$. Thus, 
\begin{align*}
\sum\limits_{\bt} \sum\limits_{b} r (\bs \cdot a, \bt \cdot b) = \;&  \sum\limits_{b} r(\bs \cdot a, \bs \cdot b) \\
= \; &   \left( p(\bs \cdot a) \wedge q(\bs \cdot a) \right) + \pi_2 (\bs \cdot a) \sum_c   d^-(\bs \cdot c) 
\sum_b \pi_1 (\bs \cdot b) \\
= \; & \left( p(\bs \cdot a) \wedge q(\bs \cdot a) \right) + d^-(\bs \cdot a) = p(\bs \cdot a) . 
\end{align*}
If $p(\bs) >  q(\bs)$, observe that $p(\bs) - q(\bs) + d^+(\bs \cdot c)  = d^-(\bs \cdot c) $ and thus
\begin{align*}
\sum\limits_{\bt} \sum\limits_{b} r (\bs \cdot a, \bt \cdot b) = \;& 
\sum\limits_{\bt \neq \bs} \sum\limits_b \pi_1 (\bt \cdot b)  \pi_2 (\bs \cdot a)  r (\bs, \bt)    +  \sum\limits_{b} r(\bs \cdot a, \bs \cdot b) \\ 
= \; & \pi_2 (\bs \cdot a) ( p(\bs) - q(\bs) ) +   \left( p(\bs \cdot a) \wedge q(\bs \cdot a) \right) 
+  \pi_2 (\bs \cdot a) \sum_c   d^+(\bs \cdot c) \\
= \; &  d^-(\bs \cdot a)  + \left( p(\bs \cdot a) \wedge q(\bs \cdot a) \right) = p(\bs \cdot a). 
\end{align*}
By symmetry, the column marginal condition holds.
\end{proof} 
 
 \begin{rem}
 Note that the above theorem continues to hold
if the common initial state $s_0$ is itself chosen randomly from
some initial probability distribution.
 \end{rem}
 
 \begin{rem}\label{rem:copula}
 MEXIT coupling is not unique in general.  
 We can (over-)parametrize all possible \MEXIT couplings by replacing the assignations  $\pi_1$ and $\pi_2$
 using copulae (\citet{Nelsen-2006}) to parametrize the dependence between changes in the \(p\)-chain and the \(q\)-chain.
 \end{rem}
 
Recall the coin flip example. 
The table for $n = 3$ given in Section 2.3 does not satisfy the inheritance principle.  
Using the construction provided in the proof above, one \MEXIT coupling is given by 

\bigskip
\begin{tabular}{|c||c|c|c|c|c|c|c|c||c|}
	\hline
	$X \backslash Y$ & HHH & HHT & HTH & HTT & THH & THT & TTH & TTT & SUM \cr
	\hline
	HHH & 0.064 & 0 & 0 & 0 & 0 & 0        &  0.0105 & 0.0505 & 0.125 \cr
	HHT & 0 & 0.096 & 0 & 0 & 0 & 0        &  0.0050 & 0.0240 & 0.125 \cr
	HTH & 0 & 0 & 0.096 & 0.019 & 0 &  0 & 0.0017 & 0.0083 & 0.125 \cr
	HTT & 0 & 0 & 0 & 0.125 & 0 & 0         & 0 & 0 & 0.125 \cr
	THH & 0 & 0 & 0 & 0  & 0.096 & 0.019  & 0.0017 & 0.0083 & 0.125 \cr
	THT & 0 & 0 & 0 & 0 & 0 & 0.125 & 0 & 0 & 0.125 \cr
	TTH & 0 & 0 & 0 & 0 & 0 & 0 & 0.125 & 0 & 0.125 \cr
	TTT & 0 & 0 & 0 & 0 & 0 & 0 & 0 & 0.125 & 0.125 \cr
	\hline
	SUM & 0.064 & 0.096 & 0.096 & 0.144 & 0.096 & 0.144 & 0.144 & 0.216 & 1 \cr
	\hline
\end{tabular}
\bigskip

It is easy to see that \MEXIT is not unique. 
Assume all the cells are fixed except the upper-right four cells, which can be seen as a $2 \times 2$ table. 
Then this $2\times 2$ table only need satisfy three constraints: the sum must be $0.9$, the sum of the first row must be $0.061$, and the sum of the first column must be $0.0155$. Hence there is still one degree of freedom.  
 
Having proven the existence of \MEXIT couplings, we now provide calculations of \MEXIT rate bounds (Subsection \ref{sec21}) and gain further insight into \MEXIT by considering its connection with the Radon-Nikodym derivative (Subsection \ref{secRN}). We finish Section \ref{sec:MEXIT} on 
an
applied note with a discussion of \MEXIT times for MCMC algorithms (Subsection \ref{secnoisy}).
\subsection{\MEXIT rate bound} \label{sec21}
We now consider \MEXIT rate bounds.

\begin{proposition} \label{prop10}
Consider the context of Theorem \ref{thm:MEXIT}. Suppose we know that
there is some $\delta>0$ such that either:\\\\
(a) for all $\bs$ and $a$,
\beqn
{ p(\bs \cdot a) / p(\bs) \over q(\bs \cdot a) / q(\bs) }
\quad \ge \quad 1-\delta
\eeqn
or \\\\
(b) for all $\bs$ and $a$,
\beqn
{ q(\bs \cdot a) / q(\bs) \over p(\bs \cdot a) / p(\bs) }
\quad \ge \quad 1-\delta.
\eeqn
Then
\beqn
\P[\text{\MEXIT at time}\,\, n+1 \, \mid \, \text{no \MEXIT by time}\,\, n] \quad \le \quad \delta.
\eeqn
\end{proposition}

\begin{proof}
Assume (a) (then (b) follows by symmetry).  We obtain
$$
\P[\hbox{\rm no \MEXIT by time } n+1 \, | \, \hbox{\rm no \MEXIT by time }n]
$$
\beqn
&=&{ \sum_{\bs,a} [p(\bs \cdot a) \wedge q(\bs \cdot a)]
\over \sum_\bs [p(\bs) \wedge q(\bs)] }\\
&\ge& 
{ \sum_{\bs,a} [(1-\delta) q(\bs \cdot a) {p(\bs) \over q(\bs)}
\wedge q(\bs \cdot a)] \over \sum_\bs [p(\bs) \wedge q(\bs)] }\\
&=&
{ \sum_{\bs,a} {q(\bs \cdot a) \over q(\bs)} [(1-\delta) p(\bs)
\wedge q(\bs)]
\over \sum_\bs [p(\bs) \wedge q(\bs)] }\\
&=&
{ \sum_\bs [(1-\delta) p(\bs) \wedge q(\bs)]
\over \sum_\bs [p(\bs) \wedge q(\bs)] }\\
&\ge& \ 1 - \delta
\, .
\eeqn
\end{proof}
The above is the discrete state-space version of a bound contained in \citet{Vollering-2016}. It should be noted that this bound applies equally well to faithful couplings, which typically degenerate in continuous time (see Theorem \ref{thm:faithful-not-maximal} in the present work for an example of this in the context of suitably regular diffusions.)
Two corollaries of Proposition \ref{prop10} follow immediately:
\begin{cor}
Under the conditions of Proposition \ref{prop10},\, $\P[$no \MEXIT by time $n] \ge (1-\delta)^n$.
\end{cor}

\begin{cor} 
Under the conditions of Proposition \ref{prop10},\, $\Expect{\text{\MEXIT time}} \ge (1/\delta).$
\end{cor}

\subsection{A Radon-Nikodym perspective on \MEXIT} \label{secRN}
\indent In this section, we explore a simple and natural connection of \MEXIT to the value of the Radon-Nikodym derivative of $q$ with respect to $p$.\\
\indent In our discussion, it will suffice to consider \MEXIT when the historical probability of the current path under both $p$ and $q$ are close to being equal, rare big jumps excepting. It follows from our \MEXIT construction that the
probability of {\it not} ``MEXITing'' by time $n$ is equal to $\sum_\bs (p(\bs) \wedge q(\bs))$, where the sum is over all
length-$n$ paths $\bs$.  Hence, conditional on having followed the path $\bs$ up to time $n$ and not ``MEXITed,'' the conditional probability of {\it not} ``MEXITing'' at time $n+1$ is equal to
$$
{ \sum_a (p(\bs \cdot a) \wedge q(\bs \cdot a))
\over
p(\bs) \wedge q(\bs) }.
\,
$$
Thus, the probability of ``MEXITing'' at time $n+1$ is
$$
1 - { \sum_a (p(\bs \cdot a) \wedge q(\bs \cdot a))
\over
p(\bs) \wedge q(\bs) }
\quad = \quad 
{ (p(\bs) \wedge q(\bs)) - \sum_a (p(\bs \cdot a) \wedge q(\bs \cdot a))
\over
p(\bs) \wedge q(\bs) }
\, .
$$
In particular, if $p(\bs) > q(\bs)$ and $p(\bs \cdot a) > q(\bs \cdot
a)$ for all $a$, then the numerator is zero, so the probability of
``MEXITing'' is zero.  That is, ``MEXITing'' can only happen when the relative
ordering of $(p(\bs), q(\bs))$ and $(p(\bs \cdot a), q(\bs \cdot a))$
are different. \\
\indent We now rephrase the above arguments in the language of Radon-Nikodym derivatives. Let
$q(a | \bs) = q(\bs \cdot a) / q(\bs)$, and
$R(\bs) = p(\bs) / q(\bs)$.  Then the non-\MEXIT probability is
$$
{ \sum_a (p(\bs \cdot a) \wedge q(\bs \cdot a))
\over
p(\bs) \wedge q(\bs) }
\quad =\quad \ \expesub{q(a|\bs)}{R(\bs \cdot a) \wedge 1 \over R(\bs) \wedge 1}
\quad = \quad \expesub{p(a|\bs)}{{R(\bs \cdot a)^{-1} \wedge 1 \over R(\bs)^{-1} \wedge 1}}
\, .
$$
Note that $\expesub{q(a|\bs)}{R(\bs \cdot a)} = R(\bs)$.  Thus, if we have either
$R(\bs)<1$ and $R(\bs \cdot a) < 1$ for all $a$, or
$R(\bs)>1$ and $R(\bs \cdot a) > 1$ for all $a$, then this non-\MEXIT
probability is one and thus the \MEXIT probability is zero.  That is, \MEXIT
can only occur when the Radon-Nikodym derivative $R$ changes from more
than 1 to less than 1 or vice-versa. 

\subsubsection{An example: \MEXIT for simple random walks} \label{sec321}
To further elucidate the connection of \MEXIT with the Radon-Nikodym derivative, we consider a concrete example: two simple random walks that both start at $0$. 
Let ``$p$'' be simple random walk with up probability $\eta < 1/2$ and down probability $1-\eta$. Similarly, let ``$q$'' be a simple random walk with
up probability $1-\eta$ and down probability $\eta$. 
The Radon-Nikodym derivative at time $n$ can be computed as 
\begin{align*}
R(\bs) = \dfrac{ p(\bs) }{q(\bs)} = \left(  \dfrac{\eta}{1-\eta} \right)^{x_n + y_n - n}  , 
\end{align*}
where $x_n$ and $y_n$ denote  the number of upward moves of chain ``$p$" and ``$q$" respectively. 
Hence $R(\bs) = 0$ if and only if $x_n + y_n = n$. Before \MEXIT, the two chains are coupled such that $x_n = y_n$, which further implies that 
\MEXIT only occurs at 0, i.e. $x_n = y_n = n/2$. %  which is the only point at which the two processes have the same transition probability, i.e.\ at which the Radon-Nikodym derivative is equal to 1. 
Indeed, the ``pre-\MEXIT'' process (i.e., the joint process, conditional on
\MEXIT not having yet occurred)
% not having ``MEXITed'') 
evolves with the following dynamics (for simplicity, we use $P$ to denote the transition probability of either chain conditional on that \MEXIT has not occurred.)

\begin{itemize}
\item For $k>0$, $P(k, k+1) = \eta$, and $P(k, k-1) = 1-\eta$.
\item For $k<0$, $P(k, k+1) = 1-\eta$, and $P(k, k-1) = \eta$.
\item $P(0, 1) = P(0, -1)=\eta$ with \MEXIT probability $1-2\eta$ when we are at $0$.
\end{itemize}
For $n = 2$, the joint distribution of the two chains is given by
\begin{center}
\begin{tabular}{|c||c|c|c|c||c|}
\hline
$q$\textbackslash$p$  &  ++  &  $+-$  &  $-+$  &  $--$ &  Sum \\
\hline 
++   &   $\eta^2$  &  0  &  0  &  $1 - 2\eta$  &  $(1 - \eta)^2$ \\  
\hline 
$+-$   &    0  &  $\eta(1-\eta)$ &  0  &  0  &  $\eta(1-\eta)$ \\
\hline
$-+$   &    0   &  0  & $\eta(1-\eta)$  & 0  &  $\eta(1-\eta)$ \\
\hline
$--$   &   0  & 0  & 0  & $\eta^2$  & $\eta^2$ \\
\hline 
Sum  &   $\eta^2$ &  $\eta(1-\eta)$ &  $\eta(1-\eta)$  & $(1-\eta)^2$ & 1 \\
\hline
\end{tabular}
\end{center}
Note that the chain $P$ is defective at 0, but otherwise has a drift towards the \MEXIT point 0. Consider the joint process, with death when MEXIT occurs. Let $Q_t$ denote the number of times this process hits 0 up to and including time $t$.  Then
\begin{equation}
\label{eq:disclocal}
\P[\text{\MEXIT by time} \,\, t \mid Q_{t-1}]\quad=\quad 1-(2\eta)^{Q_{t-1}}.
\end{equation}
Hence, $$ P[\text{no \MEXIT by time} \,\, t \mid Q_{t-1}] = (2\eta)^{Q_{t-1}}. $$ In particular, since $\eta<1/2$, and the joint process is recurrent conditional on not yet ``MEXITing'', eventual MEXIT is certain.

\subsection{An application: noisy MCMC} \label{secnoisy}
The purpose of this section is to provide an application of \MEXIT for discrete-time countable state-spaces. We do so by comparing the \MEXIT time $\tau$ of the \textit{penalty method} MCMC algorithm with the usual Metropolis-Hastings algorithm. \\
\indent In the usual Metropolis-Hastings algorithm, starting at a state $X$,
we propose a new state $Y$, and then accept it with probability
$1 \wedge A(X,Y)$, where $A(X,Y)$ is an appropriate acceptance
probability formula in order to make the resulting Markov chain reversible with respect to the {\em target} density $\pi $. In {\it noisy MCMC} (specifically, the {\it penalty method} MCMC,
see \citet{Dewing}; \citet{Fox}; \citet{Medina}; \citet{Alquier}) which is similar to but different from the {\it pseudo-marginal MCMC} method of \citet{AndrieuRoberts-2009}), we accept with probability
$\alphah(X,Y) := 1 \wedge ({\hat A}(X,Y) )$, where ${\hat A}(X,Y)$ is an
estimator of $A(X,Y)$ obtained from some auxiliary random experiment.  
\\
\indent Noisy Metropolis-Hastings is popular in situations where the target density $\pi $ is either not available or its pointwise evaluations are very computationally expensive. However replacing $A$ by ${\hat A}$ interferes with detailed balance and therefore usually the invariant distribution of noisy Metropolis-Hastings (if it even exists) is biased (ie different from $\pi $). Quantifying the bias is therefore an important theoretical question. It is not our intention to give a full analysis of this here, as this is well-studied for example \citet{Medina}.
However a crucial  component in the argument used in that paper
is the construction of a coupling between a standard and a noisy Metropolis-Hastings chain in such a way that, with high probability,  MEXIT occurs at a time after both chains have more or less converged to equilibrium. Here therefore we shall just focus on lower bounds for the MEXIT time.
\\
\indent
For this example we shall assume that $W = \exp(N)$
where $N \sim \Normal(-\sigma^2/2, \, \sigma^2)$ for some fixed $\sigma>0$
(so that $\Expect{W} = \Expect{\exp(N)} = 1$),
i.e.\ that $\alphah(X,Y) := 1 \wedge (A(X,Y) \, \exp(N))$. We now show that the \textit{penalty method MCMC} produces a Metropolis-Hastings algorithm with sub-optimal acceptance probability. 

\begin{proposition} The \textit{penalty method} MCMC
produces a Metropolis-Hastings algorithm with (sub-optimal)
acceptance probability $\alphat(X,Y, \sigma) \ :=  \Expect{\alphat(X,Y) \, | \, X,Y}$
given by
$$
\alphat(X,Y, \sigma)
\quad= \quad \Phi\left[{\log A(X,Y) \over \sigma} - {\sigma \over 2}\right]
+ A(X,Y) \, \Phi\left[-{\sigma \over 2} - {\log A(X,Y) \over \sigma} \right]
\, .
$$
\end{proposition}
\begin{proof} We invoke Proposition 2.4 of \citet{Roberts1997}, which states that if $B \sim \Normal(\mu,\sigma^2)$, then
$$
\Expect{1 \wedge e^B}
\quad = \quad \Phi\parens{{\mu \over \sigma}}+
\exp(\mu + \sigma^2/2) \, \Phi\bracks{-\sigma - {\mu \over \sigma}}
\, .
$$
Note

\beqn
\alphat(X,Y, \sigma) \quad= \quad \Expect{\alphah(X,Y)}\quad=\quad \Expect{1 \wedge (A(X,Y) e^N)}\quad=\quad \Expect{1 \wedge e^{N(-\sigma^2/2 + \log A(X,Y), \ \sigma^2)}}.
\eeqn
After straightforward algebra, the right-hand side of the last equality simplifies to

\beqn
\Phi\left[{\log A(X,Y) \over \sigma} - {\sigma \over 2}\right]
+ A(X,Y) \, \Phi\left[-{\sigma \over 2} - {\log A(X,Y) \over \sigma} \right].
\eeqn
\end{proof}

\begin{proposition} \label{prop3}
$A(X,Y) \, \phi\left[-{\sigma \over 2} - {\log A(X,Y) \over \sigma} \right]
\quad =\quad \phi\left[{\log A(X,Y) \over \sigma} - {\sigma \over 2}\right]$.
\end{proposition}

\begin{proof}
We calculate
\beqn
&&A(X,Y) \, \phi\left[-{\sigma \over 2} - {\log A(X,Y) \over \sigma} \right]\\
&=& \ {1 \over \sqrt{2\pi}} \exp\parens{ \log A(X,Y)
- {1 \over 2} \parens{-{\sigma \over 2} - \parens{\log A(X,Y) \over \sigma}^{2}}}\\
&=& \ {1 \over \sqrt{2\pi}} \exp\parens{
- {1 \over 2} \parens{ {\log A(X,Y) \over \sigma} - {\sigma \over 2} }^{2}}\\
&=&\ \phi \parens{{\log A(X,Y) \over \sigma} - {\sigma \over 2} }.
\eeqn

\end{proof}

\begin{proposition} \label{prop4}
For any $a,s > 0$, we have that
\begin{equation}
{1 \over a} \ \phi\parens{{\log a \over s} - {s \over 2}} \quad \le  \quad
{1 \over \sqrt{2\pi}}.
\end{equation}
\end{proposition}

\begin{proof}
This follows from noting
\beqn
&&{1 \over a} \ \phi\parens{{\log a \over s} - {s \over 2}} \\
&=&{1 \over \sqrt{2\pi}}
\exp\parens{- \log a -{1 \over 2} \parens{{\log a \over s} - {s \over 2}}^2} \nonumber \\
&=& {1 \over \sqrt{2\pi}} \exp\parens{-{1 \over 2} \parens{{\log a \over s} + {s \over 2}}^2}
\le {1 \over \sqrt{2\pi}} \nonumber.
\eeqn
\end{proof}

Let $r(X)$ and $\widetilde{r}(X)$ be the probabilities of rejecting the proposal when starting at $X$ for the original Metropolis-Hastings algorithm and the \textit{penalty method} MCMC, respectively. We now proceed with Proposition \ref{prop5}.

\begin{proposition} \label{prop5}
For all $X,Y$ in the state space, and $\sigma \ge 0$, the following seven statements hold\\
\hfil\break (1) $\alphat(X,Y) \le \alpha(X,Y)$.
\hfil\break (2) $\widetilde{r}(X) \ge r(X)$.
\hfil\break (3) $\lim_{\sigma \searrow 0} \alphat(X,Y, \sigma) = \alpha(X,Y)$.
\hfil\break (4) ${d \over d\sigma} \alphat(X,Y, \sigma)
= - \phi\left[{\log A(X,Y) \over \sigma} - {\sigma \over 2}\right]$.
\hfil\break (5) $0 \ge {d \over d\sigma} \alphat(X,Y, \sigma)
\ge -1/\sqrt{2\pi}$.
\hfil\break (6) $\alphat(X,Y, \sigma) \ge \alpha(X,Y) - \sigma/\sqrt{2\pi}$.
\hfil\break (7) ${\alphat(X,Y, \sigma) \over \alpha(X,Y)}
\ge 1 - \sigma/\sqrt{2\pi}$.
\end{proposition}

\begin{proof}
For statement (1), apply Jensen's inequality. Note that
\begin{align*}
\Expect{\alphat(X,Y) \, | \, X,Y}
&= \Expect{1 \wedge (A(X,Y) e^N) \, | \, X,Y}
\hspace{-1.2cm}
& \le  &
1 \wedge \Expect{(A(X,Y) e^N)} \hspace{2cm} \\
& = 1 \wedge (A(X,Y) \Expect{e^N})
& = &
1 \wedge A(X,Y) 
%\\&
\hspace{0.5cm}= \alpha(X,Y)
\, .
\end{align*}
Statement (2) follows immediately from statement (1) by taking the complements of the
expectations of the $\alpha(X,Y)$ and $\alphat(X,Y)$ with respect to $Y$.

For statement (3), note that if $A(X,Y)>1$ then
$\lim_{\sigma \searrow 0} \alphat(X,Y, \sigma) =
\Phi[+\infty] + A(X,Y) \, \Phi[-\infty]
= 1$, while if $A(X,Y)<1$ then
$\lim_{\sigma \searrow 0} \alphat(X,Y, \sigma) = 
\Phi[-\infty] + A(X,Y) \, \Phi[+\infty]
= 0  + A(X,Y) \, 1
= A(X,Y)$. Further, if $A(X,Y)=1$ then
$\lim_{\sigma \searrow 0} \alphat(X,Y, \sigma) = 
\Phi[0] + A(X,Y) \, \Phi[0]
= (1/2) + (1)(1/2)
= 1$.
Thus, in all cases,
$\lim_{\sigma \searrow 0} \alphat(X,Y, \sigma) =
1 \wedge A(X,Y)
= \alpha(X,Y)$.

For statement (4), we use Proposition \ref{prop3} to compute \\
\beqn
&& {d \over d\sigma} \alphat(X,Y, \sigma)\\
&=& {d \over d\sigma} \parens{\Phi\left[{\log A(X,Y) \over \sigma} - {\sigma \over 2}\right]
+ A(X,Y) \, \Phi\left[-{\sigma \over 2} - {\log A(X,Y) \over \sigma} \right]}\\
&=& \ \phi\left[{\log A(X,Y) \over \sigma} - {\sigma \over 2}\right]
\ \parens{- {\log A(X,Y) \over \sigma^2} - {1 \over 2}}
+ A(X,Y) \ \phi\left[-{\sigma \over 2} - {\log A(X,Y) \over \sigma} \right]\\
&=& -{1 \over 2} + {\log A(X,Y) \over \sigma^2}
\ = \ - \phi\left[{\log A(X,Y) \over \sigma} - {\sigma \over 2}\right].
\eeqn
Since $0 \le \phi(\cdot) \le {1 \over
\sqrt{2\pi}}$, statement (5) follows immediately. Statement (6) then follows by integrating from 0 to $\sigma$. For statement (7), note that if $A(X,Y) \ge 1$ then $\alpha(X,Y)=1$ and the result then follows
from statement (6).  If instead $A(X,Y) < 1$, then $\alpha(X,Y)=A(X,Y)$, and we may invoke Proposition \ref{prop4} to obtain\\
\beqn
{\alphat(X,Y,\sigma) \over \alpha(X,Y)}
\ &=& \ 1 - {\alpha(X,Y) - \alphat(X,Y,\sigma) \over \alpha(X,Y)}\\
\ &=& \ 1 - \int_{u=0}^\sigma {1 \over \alpha(X,Y)}{d \over du} \alphat(X,Y, u)du\\
&=& \ 1 - \int_{u=0}^\sigma {1 \over A(X,Y)}  \phi\left[{\log A(X,Y) \over \sigma} - {\sigma \over 2}\right] du\\
&\ge& \ 1 -  \int_{u=0}^\sigma {1 \over \sqrt{2\pi}} \ du
\quad = \quad 1 -  {\sigma \over \sqrt{2\pi}}
\, .
\eeqn
This concludes the proof.
\end{proof}
\indent Let $P$ be the law of a Metropolis-Hastings algorithm, and
$\Pt$ the law of a corresponding noisy MCMC. We now prove Proposition \ref{prop8} below, whose Corollary \ref{cor17} uses \MEXIT to control the discrepancy between the Metropolis-Hastings algorithm and the noisy MCMC algorithm.

\begin{proposition} \label{prop8}

\beqn
{d\Pt^{t+1}(\bs \cdot a) \over dP^{t+1}(\bs \cdot a)}
\  \quad \ge \quad 
{d\Pt^{t}(\bs) \over dP^{t}(\bs)}
\ \parens{1-{\sigma \over \sqrt{2\pi}}}.
\eeqn
\end{proposition}

\begin{proof}
Note first that ${d\Pt^{t}(\bs) \over dP^{t}(\bs)} = \gamma_1 \gamma_2
\ldots \gamma_n$ where each $\gamma_i$ equals either ${\alphat(X_{i-1},X_i)
\over \alpha(X_{i-1},X_i)}$ if the move from $X_{i-1}$ to
$X_i$ is accepted and otherwise ${\widetilde{r}(X) \over r(X)}$ if the move is rejected.
Statement (2) of Proposition \ref{prop5} tells us that, if we reject, 
\beqn
{d\Pt^{t+1}(\bs \cdot a) \over dP^{t+1}(\bs \cdot a)}
\quad \ge \quad 
{d\Pt^{t}(\bs) \over dP^{t}(\bs)}
\quad \ge \quad
{d\Pt^{t}(\bs) \over dP^{t}(\bs)} \parens{1-{\sigma \over \sqrt{2\pi}}}.
\eeqn
However, if we accept, then by statement (7) in Proposition \ref{prop5},
${d\Pt^{t+1}(\bs \cdot a) \over dP^{t+1}(\bs \cdot a)}
\ \ge \
{d\Pt^{t}(\bs) \over dP^{t}(\bs)} (1-{\sigma \over \sqrt{2\pi}})$, as claimed.
\end{proof}
The following Corollary to Proposition \ref{prop8} now follows immediately.
\begin{cor} \label{cor17}
${d\Pt^{t}(\bs) \over dP^{t}(\bs)}
\ \ge \ \parens{1-{\sigma \over \sqrt{2\pi}}}^t$.
\end{cor}
Applying Proposition \ref{prop8} to Proposition \ref{prop10} in Subsection \ref{sec21}, with $\delta = {\sigma \over \sqrt{2\pi}}$, the following Corollary follows immediately.

\begin{cor}
The \MEXIT time $\tau$ of the above \textit{penalty method} MCMC algorithm,
compared to the regular Metropolis-Hastings algorithm, satisfies the following two inequalities:

\beqn
\P[\tau > n]  \quad \ge  \quad \parens{1 - {\sigma \over \sqrt{2\pi}}}^n
\eeqn
and
\beqn
\Expect{\tau}  \quad \ge  \quad \sqrt{2\pi} / \sigma.
\eeqn
\end{cor}
Of course, unless $\sigma $ is small, MEXIT will likely occur substantially before Markov chain mixing, reflecting the fact that successful couplings usually have to bring chains together and not just stop them from separating. Therefore these results are usually not useful for explicitly estimating the bias of noisy Metropolis-Hastings. However they are particularly useful for demonstrating robustness results for  both noisy and pseudo-marginal chains as in
\citet{Medina} and \citet{AndrieuRoberts-2009}.

 % =========================================   
 \section{\MEXIT for general random processes}\label{sec:generalMEXIT}
The methods and results of Section \ref{sec:MEXIT} generalize to the case when the two processes are general time-inhomogeneous random processes in discrete time with countable state-space:
such processes, with state augmented to include genealogy, become Markov chains.
In fact the methods and results extend to still more general processes: in this section we deal with the case of random processes for which the state-space is a general Polish space 
(a \(\sigma\)-algebra arising from a complete separable metric space).

 % =========================================   
 \subsection{Case of one time-step}\label{sec:generalMEXIT1}
% We review the simplest case in order 
To establish notation, we first review the simplest case of just one time-step. 
We require the state-space to be Polish 
(we note that in principle one might be able to generalize a little beyond this, 
but the prospective rewards of such a generalization seem to be not very substantial).
In the case of Polish space, the diagonal set \(\Diagonal{}=\{(x,x):x\in\StateSpace\}\subset E\times E\) belongs to the product \(\sigma\)-algebra \(\SigmaAlgebra*\SigmaAlgebra\) 
(counterexamples for some more general spaces are provided in \citet[Subsection 1.6]{Stoyanov-1997}; in principle one could seek to exploit the fact that \(\Diagonal{}\) is in general
analytic with respect to \(\SigmaAlgebra*\SigmaAlgebra\), but some kind of assumption about the state-space would still be required to take care of further complications). 
\\
% However 
% \(\Diagonal{}\) 
% \emph{is} analytic and thus universally measurable.
% Consequently,
% once we have determined a coupling probability measure on the product measure space \((\StateSpace\times\StateSpace,\SigmaAlgebra*\SigmaAlgebra)\), then the diagonal set
% \(\Diagonal{}\) will belong to the completion of the product \(\sigma\)-algebra \(\SigmaAlgebra*\SigmaAlgebra\) under the coupling probability measure.
% More information about this aspect of measure-theoretic probability can be found, for example, in \cite{DellacherieMeyer-1979}.\\
\indent Consider two \(\StateSpace\)-valued random variables \(X_1^+\) and \(X_1^-\),
measurable with respect to 
% a specific sub-\(\sigma\)-algebra \(\SigmaAlgebra_1\leq\SigmaAlgebra\) 
\(\SigmaAlgebra\) on \(\StateSpace\), with distributions 
\(\Law{X_1^+}=\mu_1^+\) and \(\Law{X_1^-}=\mu_1^-\) on 
\((\StateSpace,\SigmaAlgebra)\).
% \((\StateSpace,\SigmaAlgebra_1)\).
We recall that the \emph{meet measure} \(\hat\mu_1=\mu_1^+\wedge\mu_1^-\) 
of the probability measures \(\mu^+\) and \(\mu^-\)
in the lattice of non-negative measures on \((\StateSpace,\SigmaAlgebra_1)\)
can be described explicitly using the Hahn-Jordan decomposition (\citet[\S28]{Halmos-1978}) as
\begin{equation}\label{eq:hahnjordan1}
 \mu_1^+-\mu_1^- \quad=\quad \nu_1^+-\nu_1^-
\end{equation}
for unique non-negative measures \(\nu_1^+\) and \(\nu_1^-\) of disjoint support.
The condition of disjoint support implies that
\begin{equation}\label{eq:hahnjordan2}
 \hat\mu_1 \quad=\quad \mu_1^+-\nu_1^+ \quad=\quad \mu_1^--\nu_1^-
\end{equation}
is the maximal non-negative measure \(\widetilde\mu\) such that
\[
 \widetilde\mu(D)\quad\leq\quad \min\{\mu_1^+(D), \mu_1^-(D)\} \qquad\text{ for all } D\in\SigmaAlgebra\,.
\]
% More explicitly, note that the probability measures \(\mu_1^\pm\) are both absolutely continuous with respect to their average, \(\overline\mu_1=\tfrac12(\mu_1^+ + {\mu_1^-})\). 
% This permits employment of the Radon-Nikodym theorem: the equations
% \begin{equation}\label{eq:densities}
%  f_1^\pm \quad=\quad \d\mu_1^\pm / \d\overline\mu_1
% \end{equation}
% determine non-negative densities \(f_1^\pm\) for \(\mu_1^\pm\).
% The densities are defined up to sets of \(\overline\mu_1\)-measure zero, and we suppose specific choices of non-negative representatives are made within the respective equivalence classes.
% The disjoint supports of the \(\nu_1^\pm\) can then be defined by \([f_1^\pm>1]\).
% Note also that \(f_1^++f_1^-=2\), and that \(\d\hat\mu_1=(f_1^+\wedge f_1^-)\d\overline \mu_1\),
% and that
% \begin{equation*}
%  \hat\mu_1(D) \quad=\quad \mu_1^\pm(D\cap[f_1^+<1]) + \mu_2^\pm(D\cap[f_1^+\geq1]) \qquad \text{ for all } D\in\SigmaAlgebra\,.  
%  \end{equation*}
% It is useful to note that
% \begin{equation}\label{eq:nu}
%  \d\nu_1^\pm \quad=\quad ((f_1^\pm-f_1^\mp)\vee0)\;\d\overline\mu_1\,.
% \end{equation}

\begin{lem}\label{eq:general-mexit1}
 Consider two random variables  \(X_1^+\) and \(X_1^-\) taking values in the same measurable space \((\StateSpace,\SigmaAlgebra)\) which is required to be Polish.
 The simplest \MEXIT problem 
 is solved by maximal coupling of the two marginal probability measures \( \mu_1^{+}=\Law{X_1^+}\) and \(\mu_1^-=\Law{X_1^-}\) using a joint probability measure \(m_1\) on 
the product measure space \((\StateSpace\times\StateSpace,\SigmaAlgebra*\SigmaAlgebra)\) such that
\begin{enumerate}
\item \(m_1\) has marginal distributions \(\mu_1^+\) and \(\mu_1^-\) on the two coordinates,
\item \(m_1\geq {{{\DiagonalInjection{}}_{*}}}\hat\mu_1\),
where the non-negative measure \(\hat\mu_1=\mu_1^+\wedge\mu_1^-\) is the meet measure for \(\mu_1^+\) 
and \(\mu_1^-\), and \({{{\DiagonalInjection{}}_{*}}}\) is the push-forward map corresponding to 
the \((\SigmaAlgebra:\SigmaAlgebra*\SigmaAlgebra)\)-measurable ``diagonal injection'' \({\DiagonalInjection{}}:\StateSpace\to\StateSpace\times\StateSpace\) given by \({\DiagonalInjection{}}(x)=(x,x)\).
\end{enumerate}
\end{lem}
\begin{proof}
One possible explicit construction for \(m_1\) is
 \begin{equation}\label{eq:maximal-coupling1}
  m_1 \quad=\quad {{{\DiagonalInjection{}}_{*}}}\hat\mu_1 + \frac{1}{\nu_1^+(\StateSpace)}\nu_1^+\otimes\nu_1^- \,,   
 \end{equation}
where $\nu_1^\pm$ are defined by the Hahn-Jordan decomposition in \eqref{eq:hahnjordan1} and \(\nu_1^+\otimes\nu_1^-\) is the product measure
on \((\StateSpace\times\StateSpace,\SigmaAlgebra*\SigmaAlgebra)\). 
It follows directly from \eqref{eq:hahnjordan1} that \(\nu_1^+(\StateSpace)=\nu_1^-(\StateSpace)\).
Maximality of the coupling (which is to say, maximality of \(m_1(\Diagonal{})=\hat\mu_1(\StateSpace)\)  compared to all other probability measures with these marginals)
follows from maximality of the meet measure \(\hat\mu\). This completes the proof.
\end{proof}

% \begin{comment}
% \begin{rem}\label{rem:Frechet}
%  The maximal coupling \eqref{eq:maximal-coupling1} is by no means unique:
%  the scaled product measure
%  in \eqref{eq:maximal-coupling1}
% %  \(\tfrac{1}{\nu_1^+(\StateSpace)}\nu_1^+\otimes\nu_1^-\) 
%  can be replaced by any other bivariate measure \(\Gamma_1\) in the \emph{Fr\'echet class} \(F(\nu_1^+,\nu_1^-)\)
%  of non-negative measures on \((\StateSpace\times\StateSpace,\SigmaAlgebra*\SigmaAlgebra)\) with marginals \(\nu_1^+\) and \(\nu_1^-\);
%  and the general maximal coupling can then be expressed as
%  \begin{equation}\label{eq:general-maximal-coupling}
%   m_1 \quad=\quad {{{\DiagonalInjection{}}_{*}}}\hat\mu_1 + \Gamma_1 
%  \end{equation}
%  for specified \(\Gamma_1\in F(\nu_1^+,\nu_1^-)\).
%  (Note that \(F(\nu_1^+,\nu_1^-)\) is non-void only when, as in this case, the total masses \(\nu_1^+(\StateSpace)=\nu_1^-(\StateSpace)\) agree.)
% Moreover the Fr\'echet class \(F(\nu_1^+,\nu_1^-)\) can be (over-)parametrized by the corresponding class of copulae \citep{Nelsen-2006}
% using inverse distribution transform techniques.
% \end{rem}
% \end{comment}

Given this construction, we can realize \(X_1^+\) and \(X_1^-\) as the coordinate maps for \(\StateSpace\times\StateSpace\):
% The diagonal set \(\Diagonal{}\) lies in the \(\sigma\)-algebra obtained by completing \(\SigmaAlgebra*\SigmaAlgebra\) under \(m_1\).
% Consequently 
the probability statements
\begin{equation}\label{eq:max-coupling1}
 \Prob{X_1^+\in D\;;\; X_1^+=X_1^-}\quad=\quad \hat\mu_1(D) \qquad\text{ for all } D\in\SigmaAlgebra
\end{equation}
hold for any maximal coupling of \(X_1^+\) and \(X_1^-\).

It is convenient at this point to note a quick way to recognize when a given coupling is maximal.
\begin{lem}[Recognition Lemma for Maximal Coupling]\label{lem:recognition}
Suppose the measurable space \((\StateSpace,\SigmaAlgebra)\) is Polish.
 Given a coupling probability measure \(m^*\) for \((\StateSpace,\SigmaAlgebra)\)-valued random variables \(X_1^+\) and \(X_1^-\)
 (with distributions \(\Law{X_1^+}=\mu_1^+\) and \(\Law{X_1^-}=\mu_1^-\)), 
 this coupling is maximal if
 the two non-negative measures
 \begin{equation}\label{eq:recognition1}
  \nu_1^{\pm,*} \;:\; D \quad\mapsto\quad m^*[X_1^\pm\in D\;;\; X_1^+\neq X_1^-]
 \end{equation}
(defined for \(D\in\SigmaAlgebra\)) are supported by two disjoint \(\SigmaAlgebra\)-measurable sets.
Moreover in this case the meet measure for the two probability distributions \(\Law{X_1^+}\) and \(\Law{X_1^-}\) is given by
 \begin{equation}\label{eq:recognition2}
  \hat\mu_1(D) \quad=\quad m^*[X_1^+\in D\;;\; X_1^+= X_1^-]
  \qquad \text{ for all }D\in\SigmaAlgebra\,.
 \end{equation}
\end{lem}
\begin{proof}
 This follows immediately from the uniqueness of the non-negative measures \(\nu_1^\pm\) of disjoint support appearing in the Hahn-Jordan decomposition, since
 a sample-wise cancellation of events shows that
 \[
  \mu_1^+-\mu_1^- \quad=\quad \Law{X_1^+}-\Law{X_1^-}\quad=\quad  \nu_1^{+,*} -  \nu_1^{-,*}\,.
 \]
\end{proof}

 % =========================================   
 \subsection[Case of n time-steps]{Case of $n$ time-steps}\label{sec:generalMEXITn}

 The next step is to consider the extent to which Theorem \ref{thm:MEXIT} generalizes
to the case of discrete-time random processes taking values in general Polish state-spaces.
% As noted above, \(\SigmaAlgebra*\SigmaAlgebra\) measurability of the diagonal is not required.
% Nevertheless, a
% As pointed out by an anonymous referee, 
We first note that
the generalization beyond Polish spaces cannot always hold.
Based on the work of \citet{RigoThorisson-2016}, and dating back to \citet[p.624]{Doob-1953}, 
\citet[p.210]{Halmos-1978}, and \citet[Chapter 33]{Billingsley-1968}, consider the following counterexample. \\
\indent Consider the interval $\Omega = [0, 1]$ equipped with Lebesgue measure.  
There exists a set $M \subset \Omega$ with outer measure $1$ and inner measure $0$, e.g. a Vitali set with outer measure $1$. 
Let $\mathcal{B}$ be the Borel $\sigma$-algebra on $\Omega$ and consider the $\sigma$-algebra $\sigma(\mathcal{B}, M)$. 
It can be shown that any set $A \in \sigma(\mathcal{B}, M)$ can be written as
\begin{align*}
A = (M \cap B_1 )  \cup (M^c \cap B_2), \quad B_1, B_2 \in \mathcal{B} . 
\end{align*}
The representation is not unique. 
However, using the identity $\Leb^*( M ) = \Leb^*( M \cap B_1 ) + \Leb^* ( M \cap B_1^c )$ (since $B_1$ is Lebesgue measurable), we can show  $\Leb^*( M \cap B_1) = \Leb(B_1)$
where $\Leb^*$ is the Lebesgue outer measure. 
Similarly, $\Leb^*( M^c \cap B_2) = \Leb(B_2)$. 
Hence if there is another representation $ A = (M \cap B_3 )  \cup (M^c \cap B_4)$ where $B_3$ and $B_4$ are Borel, we must have $\Leb(B_1) = \Leb(B_3)$ and $\Leb(B_2) = \Leb(B_4)$. 
Now we can define the probability measures  $m^\pm$ on $ \sigma(\mathcal{B}, M)$ by 
\begin{align*}
m^+ ( A ) =  \Leb(B_1), \quad m^- (A) = \Leb (B_2). 
\end{align*} 
It is straightforward to verify that they are probability measures.
%We claim these measures are well defined probability measures.
%To see this, let $\Leb^*$ be the Lebesgue outer measure and $\Leb_*$ be the inner measure.  
%Since $B_1$ is Borel, $\Leb^*( M ) = \Leb^*( M \cap B_1 ) + \Leb^* ( M \cap B_1^c )$.  
%But $ \Leb^* ( M \cap B_1^c ) = 1 - \Leb_* (M^c \cup B_1) \leq 1 - \Leb(B_1)$. 
%Hence $\Leb^*( M  \cap B_1  )\geq \Leb(B_1)$.  
%On the other hand, it is obvious that $\Leb^*(M \cap B_1) \leq \Leb(B_1)$.
%Thus $\Leb^*( M \cap B_1) = \Leb(B_1)$. 
Note that for any Borel set $B$, we have $m^+(B) = m^-(B) = \Leb(B)$. 
Set \(\SigmaAlgebra_1=\mathcal{B}\) and \(\SigmaAlgebra_2=\sigma(\mathcal{B},M)$. 
%(clearly $\SigmaAlgebra_1 \subset \SigmaAlgebra_2$). 
Consider two random sequences $(X_1^+, X_2^+)$ and $(X_1^-, X_2^-)$. 
Let  $X_2^\pm (\omega) = \omega$ be random variables defined on $(\Omega, \SigmaAlgebra_2, m^\pm)$. 
Let $X_1^\pm$ be defined on $(\Omega, \SigmaAlgebra_1)$ and set $X_1^\pm = X_2^\pm$
(this is allowed because the function $X(\omega) = \omega $ is Borel measurable).  
Since for any $B \in \mathcal{B}$, 
\begin{align*}
\Prob{ X_1^+ \in B} =  \Prob{ X_2^+ \in B } = m^+(B) = \Leb(B), 
\end{align*}
$X_1^\pm$ have the same law (the Lebesgue measure) and thus any realization of MEXIT would have to have 
$ \Prob{X_1^+=X_1^-}=1 $, which further implies $ \Prob{X_2^+=X_2^-}=1 $. 
On the other hand,  since $m^+(M) = 1$ and $m^-(M) = 0$, we have  $|| m^+ - m^- ||_{\rm{TV}} = 1$ w.r.t $\SigmaAlgebra_2$. 
So for any coupling of $X_2^\pm$, denoted by $(\Omega^2,  \overline{\SigmaAlgebra^2}, \bm{\mu} )$,  where $\overline{\SigmaAlgebra^2}$ denotes the completion of $\SigmaAlgebra_2 \times \SigmaAlgebra_2$ w.r.t. $\bm{\mu}$, 
we must have $\bm{\mu}(\{(\omega, \omega): \omega \in \Omega \}) = 0$. 
This gives a contradiction.

However 
% if we restrict attention to Polish state-spaces then 
the existence of MEXIT follows easily
in the case of Polish spaces,
% However the
as also noted by \citet{Vollering-2016}. 
Here follows a proof by induction.
% For the sake of completeness we summarize a variant of \citeauthor{Vollering-2016}'s argument here.
 \begin{thm}\label{thm:generalMEXIT}
  Consider two discrete-time random processes \(X^+\) and \(X^-\), begun at the same fixed initial point, 
  taking values in a measurable state-space \((\StateSpace,\SigmaAlgebra)\) which is Polish, and run up to a finite time \(n\).
  Maximal \MEXIT couplings exist.
 \end{thm}

%  \tcr{
% A general comment. I think Polish is sufficient for the existence of regular conditional probability and this is exactly what we need for the existence of MEXIT. 
% Is Polish also necessary? Perhaps not. For example, see Faden (1979) ``The existence of regular conditional probabilities". 
%  }
 
\begin{proof}
 The case \(n=1\) follows directly from the general state-space arguments of Lemma \ref{eq:general-mexit1}.
 The countable product of Polish spaces is again Polish,
%  \tcr{(I think this holds for countable product and you do need countable product)} 
 so an inductive argument completes the proof if we can establish the following.

 Suppose \(X^\pm\) are two random variables taking values in a measurable space \((\StateSpace,\SigmaAlgebra_2)\) which is Polish,
 with laws \(\mu_2^\pm\). 
 Suppose \(\SigmaAlgebra_1\subseteq\SigmaAlgebra_2\) is a sub-\(\sigma\)-algebra
 such that \((\StateSpace,\SigmaAlgebra_1)\) is also Polish,
%  \\
%  \textbf{\textcolor{red}{(Do we strictly need this Polish assumption on \((\StateSpace,\SigmaAlgebra_1)\)?!)}}
%  \tcr{(I cannot answer. It is a little bit confusing to me since Polish is the property of a topological space. Perhaps your question is whether $E$, equipped with the topology implied by $\mathcal{E}_1$, is automatically Polish. If so, I wonder is this $\sigma$-algebra $\mathcal{E}_1$ always a topology? 
% I really have no idea. 
% Perhaps it is enough to assume $E$ is Polish?)}
 and let \(\mu_1^\pm\)
 be the laws of \(X^\pm\) viewed as random variables taking values in the Polish space \((\StateSpace,\SigmaAlgebra_1)\).
 Suppose \(m_1\) is a maximal coupling with marginals \(\mu_1^\pm\)
 on \((\StateSpace\times\StateSpace,\SigmaAlgebra_1*\SigmaAlgebra_1)\).
 The claim is that there then exists a maximal coupling \(m_2\) with marginals \(\mu_2^\pm\)
 on \((\StateSpace\times\StateSpace,\SigmaAlgebra_2*\SigmaAlgebra_2)\) which equals \(m_1\) 
 when restricted to \(\SigmaAlgebra_1*\SigmaAlgebra_1\).
 
 To see this, first note from Lemma \ref{eq:general-mexit1} that \({m_1}|_{\Diagonal{}}
 ={{{\DiagonalInjection{}}_{*}}}\hat\mu_1\),
 where \(\hat\mu_1\) is the sub-probability measure given by \(\hat\mu_1=\mu_1^+\wedge\mu_1^-\).
 Moreover,
 if
 \(\hat\mu_2\) is the sub-probability measure given by
 \(\hat\mu_2=\mu_2^+\wedge\mu_2^-\),
 then we can use the infimum characterization following \eqref{eq:hahnjordan2}
 to show that \(\hat\mu_2\) satisfies \(\hat\mu_2(A)\leq\hat\mu_1(A)\) for all \(A\in\SigmaAlgebra_1\).
 Write \((1-\pi_1){\d}\hat\mu_1
 ={\d}(\hat\mu_2|_{\SigmaAlgebra_1})\) to define the \(\SigmaAlgebra_1\)-measurable random variable \(\pi_1\) (with \(0\leq \pi_1\leq 1\)) as
 \emph{the conditional probability of MEXIT immediately after time \(1\)}.
 Because \((1-\pi_1){\d}\hat\mu_1\) and \({\d}\hat\mu_2\) agree on \(\SigmaAlgebra_1\),
 and because we are working with Polish spaces, we can construct a regular conditional probability kernel
 \(\hat k_{12}(x,B)\) (a probability measure on \(\SigmaAlgebra_2\) for each fixed \(x\), and \(\SigmaAlgebra_1\)-measurable in \(x)\)
 such that 
 \begin{equation}\label{eq:hatrcp12}
  {\d}\hat\mu_2 \quad=\quad (1-\pi_1)\hat k_{12} * {\d}\hat\mu_1 \,.
 \end{equation}
%  \tcr{(Would it be helpful to clarify here that $\hat{\mu}_1$ and $\hat{\mu}_2$ are not probability measures, though it doesn't matter?)}
Similarly we can construct regular conditional probability kernels \(k^\pm_{12}(x,B)\) such that
 \begin{equation}\label{eq:rcp12}
  {\d}\mu^\pm_2 \quad=\quad k^\pm_{12} * {\d}\mu^\pm_1 \,.
 \end{equation}
 Now \((1-\pi_1){{{\DiagonalInjection{}}_{*}}}(\hat k_{12} * {\d}\hat\mu_1)={{{\DiagonalInjection{}}_{*}}}{\d}\hat\mu_2\) defines a sub-probability measure on 
 \((\StateSpace\times\StateSpace,\SigmaAlgebra_2*\SigmaAlgebra_2)\) with marginals equal to each other and given by \(\hat\mu_2\) (as a consequence of \eqref{eq:hatrcp12}).
%  defined by \tcr{(this is just~\eqref{eq:hatrcp12})}
%  \({\d}\hat\mu_2 = (1-\pi_1)\hat k_{12} * {\d}\hat\mu_1\).
 The proof of the claim will be completed if we can establish the existence of a sub-probability measure \(\Gamma_2\)
 on \((\StateSpace\times\StateSpace,\SigmaAlgebra_2*\SigmaAlgebra_2)\) with marginals defined by
 \(\mu^\pm_2-\hat\mu_2\), 
 and agreeing on \(\SigmaAlgebra_1*\SigmaAlgebra_1\) with the measure 
 defined by \({\d}m_1-(1-\pi_1){{{\DiagonalInjection{}}_{*}}}{\d}\hat\mu_1\).
 Consider 
 \[
{\d}\Gamma_2\quad=\quad(k^+_{12}\otimes k^-_{12})*({\d}m_1-(1-\pi_1){{{\DiagonalInjection{}}_{*}}}{\d}\hat\mu_1)\,,  
 \]
 where \((k^+_{12}\otimes k^-_{12})((x^+,x^-),B^+\times B^-)=k^+_{12}(x^+,B^+)\times k^-_{12}(x^-,B^-)\)
 and we use the theory of product measure to extend to a kernel of product measures
 \(k^+_{12}(x^+,\cdot)\otimes k^-_{12}(x^-,\cdot)\).
 Exactly because \((k^+_{12}\otimes k^-_{12})\) is itself a regular conditional probability kernel,
 it follows that
 \(\Gamma_2\) agrees on \(\SigmaAlgebra_1*\SigmaAlgebra_1\) with the measure 
 defined by \({\d}m_1-(1-\pi_1){\d}\hat\mu_1\).
 On the other hand, because \(\Gamma_2\) is built from appropriate product regular conditional probabilities, 
 \(\Gamma_2\) has marginals defined by 
 \(k^\pm_{12}{{\d}\mu^\pm_1 - (1-\pi_1){\d}\hat\mu}={\d}\mu_2^\pm - {\d}\hat\mu_2\) as required.
 
 In summary, the required maximal coupling at the level of \(\SigmaAlgebra_2*\SigmaAlgebra_2\)
 is defined by
 \begin{equation}
 {{{\DiagonalInjection{}}_{*}}}{\d}\hat\mu_2 + {\d}\Gamma_2
 \quad=\quad
  (1-\pi_1){{{\DiagonalInjection{}}_{*}}}(\hat k_{12} * {\d}\hat\mu_1) +
  (k^+_{12}\otimes k^-_{12})*({\d}m_1-(1-\pi_1){{{\DiagonalInjection{}}_{*}}}{\d}\hat\mu_1)\,.
 \end{equation}
%  \tcr{(I think you need the ``natural projection" argument in the previous version to justify your use of $m_1$ as a measure on $\mathcal{E}_2 \times \mathcal{E}_2$.)}
 \end{proof}

 \begin{rem}
  As in the \(n=1\) case of Lemma \ref{eq:general-mexit1},
  we can generate a whole class of maximal couplings by
  using measurable selections from Fr\'echet classes to replace
  the product regular conditional probability kernel
  \((k^+_{12}\otimes k^-_{12})*({\d}m_1-(1-\pi_1){{{\DiagonalInjection{}}_{*}}}{\d}\hat\mu_1)\).
  Equally, as in the \(n=1\) case of Lemma \ref{eq:general-mexit1},
  this clearly does not exhaust all the possibilities.
 \end{rem}

 \subsection{Unbounded and/or continuous time}\label{sec:generalMEXITextras}
  
\MEXIT for all times (with no upper bound on time) follows easily so long as the Kolmogorov Extension Theorem (\citet[\S{V.6}]{Doob-1994}) can be applied.
This is certainly the case if the state-space is Polish; we state this formally as a corollary to Theorem \ref{thm:generalMEXIT} of the previous section. 
(For an example of what can go wrong in a more general measure-theoretic context for the Kolmogorov Extension Theorem, see \citet[\S2.3]{Stoyanov-1997}.)
\begin{cor}\label{cor:generalMEXIT}
  Consider two discrete-time random processes \(X^+\) and \(X^-\), begun at the same fixed initial point, taking values in a measurable state-space \((\StateSpace,\SigmaAlgebra)\) which is Polish.
  \MEXIT couplings exist through all time.
\end{cor}

Under the requirement of Polish state-space, it is also straightforward to establish a continuous-time version of the \MEXIT result for c\`adl\`ag processes.
The result requires this preliminary elementary properties about joint laws with given marginals.

\begin{lem}\label{eq:prodlaw}
Suppose that  \(\{X_i^+\}\) and \(\{X_i^-\}\) are two collections of random variables on the probability space $(\Omega, {\cal F}, {\mathbb P})$ taking values on a metric space $(E,d)$. Suppose that  \(\{  \Law{X_i^+}  \}\) and \(\{\Law{X_i^-}\}\)  are both tight. Then  \(\{\Law{X_i^+,X_i^-}\}\) is tight on $(E\times E,{\tilde d})$
where ${\tilde d}$ denotes the Euclidean product measure $d\times d$.
\end{lem}

\begin{proof}
For any $\epsilon >0$, we can find compact sets $S^+, S^-$ such that ${\mathbb P}(X_i^+ \in S^+) > 1-\epsilon/2$ and ${\mathbb P}(X_i^- \in S^-) > 1-\epsilon/2$ for all $i$.
But $S^+ \times S^-$ is ${\tilde d}-$compact and clearly ${\mathbb P}((X_i^+, X_i^-) \in S^+\times S^-) > 1-\epsilon$, so that \(\{\Law{X_i^+,X_i^-}\}\) is tight on $(E\times E,{\tilde d})$.
\end{proof}

\begin{thm}\label{thm:continuous-MEXIT}
  Consider two continuous-time \tcr{real-valued} random processes \(X^+\) and \(X^-\), begun at the same fixed initial point, \tcr{with c\`adl\`ag paths}.
  \MEXIT couplings exist through all time.
\end{thm}
\begin{proof}
 We work first up to a fixed time \(T\).
 
 The space of c\`adl\`ag paths in a complete separable metric state-space over a fixed time interval \([0,T]\) can be considered as a Polish space (\citet[Th\'eor\`eme 1]{Maisonneuve-1972}),
 {using a slight modification of the Skorokhod metric, namely the following \emph{Maisonneuve distance}:
 if \(\tau(t):[0,T]\to[0,T]\) is a non-decreasing function determining a change of time, and if
 \(|\tau|=\sup_t|\tau(t)-t|+\sup_{s\neq t} \log\left( \tfrac{\tau(t)-\tau(s)}{t-s} \right)     \), then the Maisonneuve distance is given by
 \begin{equation}\label{eq:maisonneuve metric}
  \dist_{M}(\omega,\widetilde\omega)\quad=\quad
  \inf_\tau\{|\tau| + \dist_\StateSpace((\omega\circ\tau ) - \widetilde\omega)\}
  \,,
 \end{equation}
 where \(\omega\) and \(\widetilde\omega\) are two c\`adl\`ag paths  $[0,T] \rightarrow \tcr{\Reals}. $
\tcr{Denote this metric space, which is separable and complete, by $\mathcal{D}$.}
 
 }
%  (see the remark in \cite[\S{IV.19}]{DellacherieMeyer-1979}).%
% \NB[WSK4]{I think the fact that Polish-valued c\`adl\`ag functions form a Polish space is proved in Theorem 5.6 in Ethier and Kurtz's book \emph{Markov Processes: Characterization and Convergence}. Or look at Pollard, ``Convergence of Stochastic Processes''. Check this!}
 
 Consider a sequence of discretizations \(\sigma_n\) (\(n=1,2,\ldots\)) of time-space \([0,T]\) whose meshes tend to zero,
 each discretization being a refinement of its predecessor.
 Note that by ``discretization'' we mean an ordered sequence \(\sigma=(t_1, t_2, ...)\) where \(0<t_1<t_2<\ldots\).
  Let \(X^{\pm,n}(t)=X^\pm(\tcr{\sup}\{s\in\sigma_n:s\tcr{\leq} t\})\)  define discretized approximations of \(X^\pm\) with respect to the discretization \(\sigma_n\).
 Invoking Theorem \ref{thm:generalMEXIT}, we require \(X^{+,n}\), \(X^{-,n}\) to be maximally coupled as discrete-time random processes
 sampled only at the discretization \(\sigma_n\): since they are \tcr{constant off} \(\sigma_n\), this extends to a maximal coupling
 of \(X^{+,n}\), \(X^{-,n}\) viewed as piecewise-constant processes defined over all continuous time. 
 
 For a given c\`adl\`ag path \(\omega\), the discretization of \(\omega\) by \(\sigma_n\) converges to \(\omega\) in Maisonneuve distance.
 This follows by observing that, for each fixed \(\eps>0\), 
 the time interval \([0,T]\) can be covered by pointed open intervals \(t\in(t_-,t_+)\)
 such that \(|\omega(s)-\omega(t-)|<\eps/2\) if \(s\in(s_-,t)\) and \(|\omega(s)-\omega(t)|<\eps/2\) if \(s\in(t,s_+)\). By compactness we can select a finite sub-cover.
 For sufficiently fine discretizations \(\sigma\) we can then ensure the Maisonneuve distance between \(\omega\) and the resulting discretization is smaller than \(\eps\). Consequently, both sequences \(\{\Law{X^{+,n}}:n=1,2,\ldots\}\), \(\{\Law{X^{-,n}}:n=1,2,\ldots\}\) are tight,
 and therefore by Lemma \ref{eq:prodlaw} we know that the sequence of joint distributions \(\{\Law{X^{+,n}, X^{-,n}}:n=1,2,\ldots\}\) is also tight \tcr{in the product space $\mathcal{D}  \times \mathcal{D}$.}
% We can view \(X^{+,n}\) and \(X^+\) as constructed on the same sample space, since one is the discretization of the other.
% Consequently it follows from this construction that \(X^{+,n}\to X^+\) almost surely and therefore \(X^{+,n}\Rightarrow X^+\) weakly and
% \(X^{-,n}\Rightarrow X^-\) weakly. 
% Accordingly, 

%the sequences of marginal distributions \(\{X^{+,n}:n=1,2,\ldots\}\), \(\{X^{-,n}:n=1,2,\ldots\}\) are both tight,
% and therefore 
 Therefore (selecting a weakly convergent subsequence if necessary) we may suppose the joint distribution $(X^{+,n}, X^{-,n})$ converges weakly 
\tcr{in $\mathcal{D}  \times \mathcal{D}$}
to a limit which we denote by $(\tilde{X}^+, \tilde{X}^-)$.
\tcr{
Since $(X^{+,n}, X^{-,n})$ has been constructed to satisfy \MEXIT for \(t\in \sigma_n\), and since $(X^{+,n}, X^{-,n})$ is constant off \(\sigma_n\), it follows for all \(t\) that 
\begin{equation*}%\label{eq:discrete-mexit}
\begin{aligned}
& \Prob{X^{+,n}(s) = X^{-,n}(s) \text{ for all }s<t} \\
=\quad &  \left(\Law{(X^{+,n}(s):s< t)}\wedge\Law{(X^{-,n}(s):s< t)}\right)(\Reals)  \quad=\quad  m_n(t)\,.
\end{aligned}
\end{equation*}
Let $ m_\infty(t)$ be defined analogously for \(\tilde{X}^+\) and \(\tilde{X}^-\) and note that \(m_n(t)\), \(m_\infty(t)\)
are both decreasing in \(t\); 
moreover
\begin{equation*}%\label{eq:bound}
m_n(t) \quad\downarrow \quad m_\infty(t) \qquad \text{ for } t\in \bigcup_m \sigma_m\,, 
\end{equation*}
since the left-hand side corresponds to the less onerous ``MEXIT on \(\sigma_n\)'' requirement that  
\(X^{+, n}\) and \(X^{-, n}\) be constructed to agree only on \(\sigma_n\cap[0,t)\) (a set of time points increasing in \(n\)) rather than all of \([0,t)\). 
We require the discretizations \(\sigma_n\) to be augmented (modifying $(X^{+,n}, X^{-,n})$ accordingly) so that the decreasing function \(m_\infty\) is continuous off $\cup_n \sigma_n$.
\\
\indent We now make a key observation: \MEXIT questions can be re-expressed in terms of continuous sample-path processes rather than \cadlag processes.
For \(\eps>0\), consider the smoothing operator \(S_\eps\) acting on \(f\in\mathcal{D}\) as follows
\[
 S_\eps(f)(t) \quad=\quad \frac{1}{\eps} \int_{t-\eps}^t f(u) \d{u}\,,
\]
 where we take \(f(t)=f(0)\) for \(t\leq0\). Then \(S_\eps:\mathcal{D}\to C([0,1])\) is continuous,
where  \(C([0,1])\) is the space of continuous real-valued functions on \([0,1]\),
endowed with the supremum metric.
  It therefore follows that \tcr{in $C([0,1])  \times C([0,1])=C([0,1])^2$ we have}
\[
 \left(S_\eps( X^{+,n}),S_\eps(  X^{-,n})\right) \quad\Rightarrow\quad
 \left(S_\eps(\tilde{X}^+),S_\eps(\tilde{X}^-)\right)\,.
\]
% when viewed as probability measures in \(C([0,1]\to\Reals^2)\). \\
\indent 
On the other hand, for any \(t\in[0,1]\) it follows by construction and the \cadlag property of \(f\) and \(g\) that \(S_\eps(f)(s)=S_\eps(g)(s)\) for all \(s\leq t\)
if and only if \(f(s)=g(s)\) for all \(s<t\).
} Suppose time $t$ belongs to one of the discretizations in the sub-sequence, and thus eventually to all (since each discretization is a refinement of its predecessor). 
Consider the subspace of \tcr{\(\mathcal{D}\times\mathcal{D}\)} given by $A_t = [ \mathsf{MEXIT} \geq t ]$.
\tcr{Since \([S_\eps(X^{+,n})(s)=S_\eps( X^{-,n})(s) \text{ for }s\leq t]\) and
\([S_\eps(\tilde{X}^+)(s)=S_\eps(\tilde{X}^-)(s) \text{ for }s\leq t]\)
can be viewed as corresponding to the same closed subset of \(C([0,1])^2\),}
by the Portmanteau Theorem of weak convergence (Billingsley [4, Theorem 2.1]),  
\begin{equation*}
\limsup\limits_{n \rightarrow \infty} \Prob{ (X^{+, n}, X^{-, n}) \in A_t } 
\leq \Prob{ (\tilde{X}^{+}, \tilde{X}^{-}) \in A_t }. 
\end{equation*}
\tcr{
Considerations of total variation distance tell us that  
\( \mathbb{P}[  (\tilde{X}^{+}, \tilde{X}^{-}) \in A_t ] \leq m_\infty(t)\);
indeed \(\tilde{X}^+\) and \(\tilde{X}^-\)
cannot disagree at a slower rate than that afforded by \MEXIT.
On the other hand, \(\Prob{ (\tilde{X}^{+}, \tilde{X}^{-}) \in A_t }\) relates to total variation distance as above,
% we can apply \eqref{eq:discrete-mexit},
% and so deduce
so
\begin{equation*}
 \limsup_{n\to\infty}m_n(t) \quad\leq\quad  
 \Prob{ (\tilde{X}^{+}, \tilde{X}^{-}) \in A_t } \quad\leq\quad
 m_\infty(t)
\qquad\text{ for all }t\,.
\end{equation*}
But \(m_n\downarrow m_\infty\) on \(\sigma_m\), so
% Using \eqref{eq:bound} it follows that 
\( \mathbb{P}[  (\tilde{X}^{+}, \tilde{X}^{-}) \in A_t ] = m_\infty(t)\) 
for all $t \in \cup_n \sigma_n$. 
}
The \cadlag property \tcr{and the continuity of $m_\infty$ off $\cup_n \sigma_n$} then implies maximality of the limiting coupling for all times $t \leq T$. 
\tcr{Hence $(\tilde{X}^+, \tilde{X}^-)$ is a \MEXIT construction as required. 
\MEXIT for all time follows using the Kolmogorov Extension Theorem as above. }
\end{proof}

 \begin{rem}
  \citet{SverchkovSmirnov-1990} prove a similar result for maximal couplings by means of general martingale theory.
 \end{rem}
 
 \begin{rem}
  Note that Th\'eor\`eme 1 of \citet{Maisonneuve-1972} can be viewed as justifying the notion of the space of c\`adl\`ag paths: 
  this space is the completion of the space of step functions under the Maisonneuve distance \(\dist_M\).
  Thus in some sense Theorem \ref{thm:continuous-MEXIT} is a maximally practical result concerning \MEXIT!
 \end{rem}

% 
% % =========================================   
%  \section{\MEXIT for coin-flipping}\label{sec:coin}
%   \TBC{Develop section 3 of notes by Ernst et al.}
  
% =========================================   
 \section{\MEXIT for diffusions}\label{sec:diffusions}
 The results of Section \ref{sec:generalMEXIT} apply directly to diffusions, which therefore exhibit \MEXIT. This section discusses the solution of a \MEXIT problem for Brownian motions,
which can be viewed as the limiting case for random walk \MEXIT problems.
%  \xNB[WSK2]{Need to explain why a -- non-genealogical -- \MEXIT must be maximal if the post-\MEXIT processes are confined to disjoint regions of state-space. I'm a bit puzzled about the r\^ole of genealogies here!}%
%  \xNB[WSK2]{The general argument about \MEXIT for semimartingale diffusions should be adapted to general \cadlag processes, and moved to a separate small section following Section \ref{sec:MEXIT}.}

 It is straightforward to show that \MEXIT will generally have to involve constructions not adapted to the shared filtration of the two diffusion in question.
 By ``faithful'' \MEXIT we mean a \MEXIT construction which generates a coupling between the diffusions which is Markovian with respect to the joint and individual filtrations (see \citet{Rosenthal-1997} and \citet{Kendall-2013a} for further background).
 We consider the case of elliptic diffusions \(X^+\) and \(X^-\) with continuous coefficients.
 \begin{thm}\label{thm:faithful-not-maximal}
  Suppose \(X^+\) and \(X^-\) are coupled elliptic diffusions, thus with continuous semimartingale characteristics given by their drift vectors and volatility (infinitesimal quadratic variation) matrices, begun at the same point, 
  with this initial point lying in the open set where either or both of the drift and volatility characteristics disagree. 
  Faithful \MEXIT must happen immediately.
 \end{thm}
\begin{proof}
   Let \(T\) be the \MEXIT time, which by faithfulness will be a stopping time with respect to the common filtration.
  If \(X^+\) and \(X^-\) are semimartingales agreeing up to the random time \(T\), 
        then the localization theorems of stochastic calculus tell us that the integrated drifts and quadratic variations of \(X^+\) and \(X^-\)
        must also agree up to time \(T\).
%         (I believe this result occurs in Meyer's famous Strasbourg lecture notes on stochastic calculus.)
        It follows that \(X^+\) and \(X^-\) agree as diffusions up to time \(T\).
        Were the faithful \MEXIT stopping time to have positive chance of being positive then the diffusions would have to agree on the range of the common diffusion up to faithful \MEXIT;
        this would contradict our assertion that the initial point lies in the open set where either or both of the drift and volatility characteristics disagree.
%         (In case of hypoelliptic diffusions, this implies the diffusions must agree on a finely open set including the common starting point -- ``finely'' with respect to the common pre-\MEXIT diffusion.
\end{proof}

 By way of contrast, \MEXIT can be described explicitly in the case of two real Brownian motions \(X^+\) and \(X^-\)  with constant but differing drifts.
 Because of re-scaling arguments in time and space, there is no loss of generality in supposing that both \(X^+\) and \(X^-\) begin at \(0\), with \(X^+\) having drift \(+1\) and \(X^-\) having drift \(-1\).
 \begin{thm}\label{thm:BM-MEXIT}
  If \(X^\pm\) is Brownian motion begun at \(0\) with drift \(\pm1\), then \MEXIT between \(X^+\) and \(X^-\) exists and is almost surely positive.
 \end{thm}

\begin{proof}%[Proof (outline)]
The existence of \MEXIT directly follows from Theorem 27. 
	The almost surely positiveness will be shown in Subsection \ref{sec:explicit} below, 
through a limiting version of
the 
random walk
argument in Subsection \ref{sec321}. 
Alternatively one can argue succinctly and directly
% in the case of Brownian motions \(X^\pm\) of drifts \(\pm1\),
using
% one can employ 
the excursion-theoretic arguments of \citeauthor{Williams-1974}' \citeyearpar{Williams-1974}
celebrated path-decomposition of Brownian motion with constant drift
(an exposition in book form is given in \citet{RogersWilliams-1987}).

%         It follows from this work,
%         \textcolor{blue}{and \cite{KaratzasShreve-1987}'s identification of Brownian motion time-changed only to be positive (respectively negative) to be reflected Brownian motion
%         (respectively the negative of Brownian motion),}
%         that 
        
Calculation
% A straightforward calculation 
shows that 
the \emph{bounded} positive excursions of \(X^+\) (respectively \(-X^-\)) from \(0\) are those of 
the positive excursions of a Brownian motion of \emph{negative} drift \(-1\),
while the \emph{bounded} negative excursions of \(X^+\) (respectively \(-X^-\)) from \(0\) are those of 
the negative excursions of a Brownian motion of \emph{positive} drift \(+1\).
(The unbounded excursion of \(X^+\) follows the law of the distance from its starting point of Brownian motion in hyperbolic \(3\)-space,
while the unbounded excursion of \(X^-\) has the distribution of the mirror image of the unbounded excursion of \(X^+\).)

Viewing \(X^\pm\) as generated by Poisson point processes of excursions indexed by local time,
it follows that we may couple \(X^+\) and \(X^-\) 
% in the same way as we coupled \(Y^+\) and \(Y^-\).
to share the same bounded excursions, with unbounded excursions being the reflection of each other in \(0\).
Moreover 
% once more
the processes have disjoint support once they become different. So
the Recognition Lemma for Maximal Coupling (Lemma \ref{lem:recognition}) applies,
and
hence this 
% too 
is a \MEXIT coupling.
%
%         Consequently, the processes \(X^+\) and \(X^-\) must agree in law up to the start of their unbounded excursions; moreover the \(X^-\) unbounded excursion is the reflection in law of the \(X^+\) unbounded excursion.
%         Hence we may construct a \MEXIT coupling in which \(X^+=X^-\) up to a (non-Markovian) time \(M\) when \(X^+_M=X^-_M=0\), at which point both commence their unbounded excursions and one unbounded excursion  is the reflection of the other.
%         
%         Finally, this exit coupling must indeed be maximal: the supports of the post-\MEXIT processes are disjoint, and therefore Lemma \ref{lem:recognition} applies.
\end{proof}

\subsection{Explicit calculations for Brownian \MEXIT}

%As before, 
Let \(X^+\) and \(X^-\) begin at \(0\), with \(X^+\) having drift \(+\theta \) and \(X^-\) having drift \(-\theta \) with $\theta > 0$. The purpose of this section is to offer explicit calculations of 
% non-maximal \MEXIT means and 
\MEXIT and \MEXIT means. \\
\noindent \textbf{Calculation 1.} \label{calc1}
The meet of the distributions of \(X^+_t\) and \(X^-_t\) is the meet of \(N(\theta t,t)\) and \(N(- \theta t,t)\),  and the probability of \MEXIT happening after time \(t\) is given by the total mass of this meet sub-probability distribution. Therefore:

\beqn
\Prob{\MEXIT \geq t}\quad &=&\quad \Prob{N(0,t)< - \theta t}+\Prob{N(0,t) > \theta t}\\
\quad &=& \quad 2\Prob{N(0,t)> \theta t}\\
\quad &=& \quad \frac{2}{\sqrt{2 \pi}} \int_{\theta \sqrt{t}}^\infty e^{-u^2/2} du.
\eeqn
Thus, 
\beqn
\Expect{\MEXIT}\quad =\quad \frac{2}{\sqrt{2\pi}}\int_0^\infty \int_{\theta \sqrt{t}}^\infty e^{-u^2/2} du dt= \theta^{-2}.
\eeqn

\begin{rem}
Excursion theoretic arguments can be used to confirm this is mean time to \MEXIT for the specific construction given in Theorem \ref{thm:BM-MEXIT}.
\end{rem}

% \noindent \textbf{Calculation 3.} \label{calc3}
\noindent \textbf{Calculation 2.} \label{calc2}
We now consider the expected amount of time $T$ during which processes agree before \MEXIT happens.
%We now generalize Calculation 1 by considering the maximal coupling length $L$ for general diffusions $X_t = W_t + \theta t$ and $Y_t = W_t - \theta t$. 
\begin{align*}
\Expect{T} \quad &= \quad  \int_0^\infty \expesub{W}{\min\{ e^{ \theta W_t - \theta^2 t/2}, e^{- \theta W_t - \theta^2 t/2} \}}dt  \\
\quad &=\quad  2 \int_0^\infty  \expesub{W}{ e^{- \theta W_t - \theta^2 t/2} ; W_t > 0} dt  \\
\quad &= \quad 2\int_0^\infty \int_0^\infty  \dfrac{1}{\sqrt{2\pi t}} \exp \parens{- \dfrac{ (w + \theta t)^2 }{2t}} dw dt   \\
\quad &=\quad  \theta^{-2}.
\end{align*}

\subsection{An explicit construction for \MEXIT for Brownian motions with drift}\label{sec:explicit}
\def\Proba{\mathbb{P}}

In this section, we continue the scenario of 
Calculation 2
% Calculation 3
above. We see that \MEXIT should have the complementary cumulative distribution function
\begin{equation}
\label{eq:BMdMEXIT}
\Prob{\MEXIT \geq t} = 2\Phi (-\theta \sqrt{t} ),
\end{equation}
where $\Phi (y) = \int_{-\infty }^y (2\pi )^{-1/2} e^{-u^2/2} du$.
A natural question to ask is as follows: how can one explicitly construct and understand this \MEXIT time in a way that relates to the random walk constructions of Subsection \ref{sec321}? 
In this section we first deduce a candidate coupling and EXIT time, and then we proceed to show by direct calculation that our construction indeed gives the correct \MEXIT time distribution above.

We note from the discrete constructions of Section 
\ref{sec:MEXIT} (in particular Subsection \ref{secRN}) that \MEXIT is only possible when the Radon-Nikodym derivative between the ``$p$'' and ``$q$'' processes moves from being below $1$ to above $1$ or moves from being above $1$ to below $1$. Let $\Proba ^+$, $\Proba^-$ denote the probability laws of $X^+$, $X^-$ respectively. We  have that
$$
{d\Proba^+ \over d\Proba^-}(W_{[0,T]}) =\exp \{ 2 \theta W_T    \},
$$
which is continuous in time with probability $1$ under both $\Proba ^+$ 
and $\Proba ^-$. By analogy to the discrete case, the region in which \MEXIT could possibly occur corresponds to the interface $
{d\Proba^+ \over d\Proba^-}(W_{[0,T]}) =1$ (that is, where $W_T=0$).

Now we shall focus on the random walk example at the end of Subsection \ref{secRN}. We note that the \MEXIT distribution given in (\ref{eq:disclocal})
can be constructed as the first time the occupation time of $0$ exceeds a geometric random variable with ``success'' probability $1-2\eta $. We aim to give a similar interpretation for the Brownian motion  case. To do this, we shall use a sequence of random walks converging to the appropriate Brownian motions. To this end, let
$$
\eta_n ={1\over 2}\left(1-{\theta \over n} \right),
$$
and set $X^{n+}$ and $X^{n-}$ to be the respective simple random walks with up probability $1-\eta _n$ and $\eta _n$ and sped up by factor $n^2$. We assume (unless otherwise stated) that all processes begin at $0$ so that we have that 
$$
X^{n+}(t) = \sum_{i=1}^{\lfloor n^2 t\rfloor } Z_i^{n+},
$$
where $\{X^{n+}_i\}$ denote dichotomous random variable taking the value $+1$ with probability $1-\eta _n$ and $-1$ with probability $\eta _n$. We define $X^{n-}$ analogously.

Given this setup, we have the classical weak convergence results that the law of
$X^{n+}$ converges weakly to that of $X^+$, and similarly $X^{n-}$ converges weakly to $X^-$. Moreover the joint pre-\MEXIT process described in Subsection \ref{secRN} will have drift $-sgn(X_t) \theta $. The following holds for the \MEXIT probability in
(\ref{eq:disclocal})
$$
\Prob{MEXIT >t} = \left( 1-{\theta \over n} \right)^{n \ell^{n}_t} \longrightarrow
e^{-\theta \ell ^n_t},
$$
where $\ell^{n}_t$ is the Local Time at $0$ of the pre-\MEXIT process for the $n$th  approximation random walk.

In the (formal) limit as $n \to \infty $, this 
recovers the
% gives a 
construction 
in Theorem \ref{thm:BM-MEXIT}
of
% for the 
Brownian motion \MEXIT time,
as follows. Let $X$ be the diffusion with drift $-\operatorname{sgn}(X) \theta $ and unit diffusion coefficient started at $0$ and let $\ell_t$ denote its local time at level $0$ and time $t$. Then set $E$ to be an exponential random variable with mean $\theta ^{-1}$. Then the pre-\MEXIT dynamics are described by $X$ until
$\ell_t>E$ at which time \MEXIT occurs. $E > 0$ w.p. 1 and hence \MEXIT is positive a.s. since the local time is a continuous process.

We shall now verify that this construction does indeed achieve the valid \MEXIT probability given in (\ref{eq:BMdMEXIT}).  By integrating out $E$ we are required to show that
$$
\Expect{e^{-\theta \ell _t}} = 2 \Phi (-\theta \sqrt{t} )\ .
$$

We proceed to do so. Firstly, we note that by symmetry, we may set $\ell_t$ to be the local time at level $0$ of Brownian motion with drift $-\theta $ reflected at $0$. Note that by an extension of L\'evy's Theorem (see \citet{peskir2006reflecting}) that the law of $\ell_t$ is the same as that of the maximum of Brownian motion with drift $\theta $, i.e. that of
$X^+$. Now this law is well-known as the Bachelier-L\'evy formula (see for example \citet{Lerche2013boundary}):
$$
\Prob{\ell_t<a} = \Phi \left({a \over \sqrt{t}}-\theta \sqrt{t}    \right) 
- e^{2 a \theta }
\Phi\left({-a \over \sqrt{t}}-\theta \sqrt{t}    \right),
$$
with density
$$
f_{\ell_t}(a)={1\over \sqrt{t}} \left(
 \phi \left({a \over \sqrt{t}}-\theta \sqrt{t}    \right) + e^{2 a \theta }
\phi\left({-a \over \sqrt{t}}-\theta \sqrt{t}    \right) -2\sqrt{t}\theta e^{2 a \theta }
\Phi\left({-a \over \sqrt{t}}-\theta \sqrt{t}    \right)
\right),
$$
where $\phi $ is the standard normal density function $\phi (y) = (2\pi ) ^{-1/2} e^{-y^2/2}$. By direct manipulation of the exponential quadratic in the second of the three terms above, it can readily be shown to equal the first term. Thus
$$
f_{\ell_t}(a)={1\over \sqrt{t}} \left(
 2 \phi \left({a \over \sqrt{t}}-\theta \sqrt{t}    \right) 
  -2\sqrt{t}\theta e^{2 a \theta }
\Phi\left({-a \over \sqrt{t}}-\theta \sqrt{t}    \right)
\right).
$$
We now directly calculate the Laplace transform of this distribution to obtain
(\ref{eq:BMdMEXIT}).
$$
\Expect{
e^{-\theta \ell _t}
} ={2 \over \sqrt{t} } \int_0^{\infty } e^{-\theta a} \left(
\phi \left({a \over \sqrt{t}}-\theta \sqrt{t}    \right) 
  -\sqrt{t} \theta e^{2 a \theta }
\Phi\left({-a \over \sqrt{t}}-\theta \sqrt{t}    \right) \right) da
$$
$$
=:  {2 \over \sqrt{t} } (T_1 - T_2)\  .
$$
Using integration by parts, we easily work with $T_2$ to obtain
$$
T_2 = T_1 -  \sqrt{t} \Phi (-\theta \sqrt{t} ),
$$
which implies the assertion in (\ref{eq:BMdMEXIT}), as required.

%\begin{rem}
 %A similar construction to that in Remark \ref{remmeet} works for simple asymmetric random walk and its opposite under reflection, %and the construction can be analysed using a similar analysis,
%albeit without the benefit of Brownian scaling.
%\end{rem}

 % =========================================    
 \section{Conclusion}\label{sec:conclusion}
In this paper, we have studied an alternative coupling framework in which one seeks to arrange for two different Markov processes to remain equal for as long as possible, when started in the same state. We call  this ``un-coupling'' or ``maximal agreement'' construction \emph{\MEXIT}, standing for ``maximal exit'' time. \MEXIT sharply differs from the more traditional maximal coupling constructions 
studied in \citet{Griffeath-1975}, \citet{Pitman-1976}, and \citet{Goldstein-1978} in which one seeks to build two different
copies of the same Markov process started at two different
initial states in such a way that they become equal as soon as possible.

\indent This work begins with practical motivation for \MEXIT by highlighting the importance of un-coupling/maximal agreement arguments 
in a few key statistical and probabilistic settings. With this motivation established,
we develop an explicit \MEXIT construction for Markov chains in discrete time with countable state-space.
We then generalize the construction of \MEXIT to random processes on Polish state-space in continuous time.
We conclude with the solution of a \MEXIT problem for Brownian motions.\\
\indent As noted in Remark \ref{rem:copula}, the approach that we have followed in the construction of 
\MEXIT introduces the role of copula theory in parametrising varieties of maximal couplings for random processes. 
Our future work will aim to establish a definitive role for \MEXIT (as well as for probabilistic coupling theory in general) in copula theory.\\

\noindent \textbf{Acknowledgments}\\
% We thank the anonymous referees for their valuable suggestions
% which greatly improved the quality of this manuscript. 
We are grateful to the two anonymous referees whose suggestions helped improve the quality of this manuscript. The first author thanks Professor Larry Shepp for talking with him in 2011-2012 about the ideas that would eventually become MEXIT. He also thanks Dr. Quan Zhou for kind discussion.
The second author thanks Mateusz Majka and Aleks Mijatovi\'c for useful conversations, and acknowledges support by EPSRC grant EP/K013939 and also by the Alan Turing Institute under the EPSRC
grant EP/N510129/1. 
The third author was supported by EPSRC grants EP/K034154/1, EP/K014463/1, EP/D002060/1. 
The fourth author received research support from NSERC of Canada.

\end{document}